\documentclass[reqno,12pt]{amsart}

\usepackage{amsmath,amsthm,amssymb,amsfonts, mathabx, tikz-cd}
\usepackage{commath}
\usepackage{extarrows}

\usepackage{tikz}
\usetikzlibrary{decorations.pathreplacing}
\usepackage[dvipsnames]{xcolor}
\usepackage{tabularx}
\usepackage{longtable}
\usepackage{array}
\usepackage{faktor}

\usepackage{catex}
\usepackage[all]{xy}

\usepackage{enumitem}
\usepackage{float}

\usepackage{titletoc}
\usepackage[titletoc]{appendix}

\usepackage{geometry}

\usepackage[hyperfootnotes=false,draft=false,pdftex, colorlinks=true, linkcolor=black, citecolor=black, urlcolor=black, hypertexnames=false]{hyperref}

\usepackage{adjustbox}

\usepackage{combelow}

\numberwithin{equation}{section}

\newcommand{\eq}{\,=\,}

\newcommand{\ms}[1]{\mathsf{#1}}

\newcommand\Simeq{\;\simeq\;}

\newcommand{\overbar}[1]{\mkern 1mu\overline{\mkern-1mu#1\mkern-1mu}\mkern 1mu}
\newcommand{\highoverline}[1]{\overbar{\vphantom{#1^{\rule{0pt}{0.3ex}}}#1}}

\newcommand{\rotamp}[2]{\mathbin{\rotatebox[origin=c]{180}{$#1\&$}}}
\newcommand{\parll}{{\mathpalette\rotamp\relax}}
\DeclareMathOperator{\partimes}{{\parll}}
\newcommand{\ot}{\otimes}

\newcommand{\tim}{\times}

\newcommand{\id}{{\operatorname{id}}}
\newcommand{\idC}{{\operatorname{id}_{\cC}}}

\newcommand{\ao}{\alpha}

\newcommand{\lo}{\lambda}

\newcommand{\ro}{\rho}

\newcommand{\aoc}{a}

\newcommand{\loc}{l}

\newcommand{\roc}{r}

\newcommand{\ap}{\highoverline{\alpha}}

\newcommand{\lp}{{\highoverline{\lo}}}

\newcommand{\rp}{{\overline{\ro}}}

\newcommand{\apc}{\overline{a}}

\newcommand{\lpc}{{\overline{l}}}

\newcommand{\rpc}{{\overline{\roc}}}

\newcommand{\distl}{\delta^l}
\newcommand{\distr}{\delta^r}

\newcommand{\distlc}{d^l}
\newcommand{\distrc}{d^r}

\newcommand{\slongrightarrow}[1]{\,\longrightarrow{#1}\,}
\newcommand{\sxlongrightarrow}[1]{\,\xlongrightarrow{#1}\,}

\newcommand{\oLD}{\highoverline{\ms{LD}}}
\newcommand{\Ord}{\mathbf{Ord}}

\newcommand{\Neq}{\mathrel{\scalebox{1.1}{$\neq$}}}
\newcommand{\NNez}{\NN_{0}}

\def\cC{\mathcal C}

\def\CC{\mathbb C}

\def\NN{\mathbb N}

\def\ZZ{\mathbb Z}
\def\Si{\Sigma}
\def\wSi{\widehat{\Sigma}}

\theoremstyle{plain}
\newtheorem{theorem}{Theorem}[section]
\newtheorem*{maintheorem*}{Main Theorem}
\newtheorem{lemma}[theorem]{Lemma}
\newtheorem{proposition}[theorem]{Proposition}

\theoremstyle{definition} 
\newtheorem{definition}[theorem]{Definition}

\theoremstyle{remark}
\newtheorem{remark}[theorem]{Remark}
\newtheorem{example}[theorem]{Example}

\usepackage{scalerel}
\newcommand{\mrhup}{\scaleobj{0.5}{\rightharpoonup}}
\newcommand{\mlhup}{\scaleobj{0.5}{\leftharpoonup}}
\newcommand{\mrelb}{\hstretch{0.5}{\scaleobj{0.5}{\relbar}}}
\makeatletter
\def\rhfill@#1{$\m@th\thickmuskip0mu\medmuskip\thickmuskip\thinmuskip\thickmuskip
	\relax#1\mkern+2.5mu\cleaders\hbox{$#1\mrelb\mkern-1mu$}\hfill\mkern-6mu\mrhup\mkern+0.25mu$}
\def\lhfill@#1{$\m@th\thickmuskip0mu\medmuskip\thickmuskip\thinmuskip\thickmuskip
	\relax#1\mkern+1.5mu\mlhup \cleaders\hbox{$#1\mkern-2.5mu\mrelb$}\hfil\mkern+0.5mu$}
\DeclareRobustCommand{\overrightharpoon}{\mathpalette{\overarrow@\rhfill@}}
\DeclareRobustCommand{\overleftharpoon}{\mathpalette{\overarrow@\lhfill@}}
\makeatother

\let\prt=\partimes

\usepackage[hyperfootnotes=false,draft=false,pdftex, colorlinks=true, linkcolor=black, citecolor=black, urlcolor=black, hypertexnames=false
]{hyperref}

\UseTips
\SelectTips{eu}{11}

\newdir{ >}{{}*!/-10pt/@{>}}
\newdir{ -}{{}*!/-10pt/@{}}
\newdir{> }{{}*!/+10pt/@{>}}

\makeatletter

\xyletcsnamecsname@{dir4{}}{dir{}}
\xydefcsname@{dir4{-}}{\line@ \quadruple@\xydashh@}
\xydefcsname@{dir4{.}}{\point@ \quadruple@\xydashh@}
\xydefcsname@{dir4{~}}{\squiggle@ \quadruple@\xybsqlh@}
\xydefcsname@{dir4{>}}{\Tttip@}
\xydefcsname@{dir4{<}}{\reverseDirection@\Tttip@}

\xydef@\quadruple@#1{%
	\edef\Drop@@{%
		\dimen@=#1\relax
		\dimen@=.5\dimen@
		\A@=-\sinDirection\dimen@
		\B@=\cosDirection\dimen@
		\setboxz@h{%
			\setbox2=\hbox{\kern3\A@\raise3\B@\copy\z@}%
			\dp2=\z@ \ht2=\z@ \wd2=\z@ \box2
			\setbox2=\hbox{\kern\A@\raise\B@\copy\z@}%
			\dp2=\z@ \ht2=\z@ \wd2=\z@ \box2
			\setbox2=\hbox{\kern-\A@\raise-\B@\copy\z@}%
			\dp2=\z@ \ht2=\z@ \wd2=\z@ \box2
			\setbox2=\hbox{\kern-3\A@\raise-3\B@ \noexpand\boxz@}%
			\dp2=\z@ \ht2=\z@ \wd2=\z@ \box2
		}%
		\ht\z@=\z@ \dp\z@=\z@ \wd\z@=\z@ \noexpand\styledboxz@
	}%
}

\xydef@\Tttip@{\kern2pt \vrule height2pt depth2pt width\z@
	\Tttip@@ \kern2pt \egroup
	\U@c=0pt \D@c=0pt \L@c=0pt \R@c=0pt \Edge@c={\circleEdge}%
	\def\Leftness@{.5}\def\Upness@{.5}%
	\def\Drop@@{\styledboxz@}\def\Connect@@{\straight@{\dottedSpread@\jot}}}

\xydef@\Tttip@@{%
	\dimen@=.25\dimen@
	\B@=\cosDirection\dimen@
	\setboxz@h\bgroup\reverseDirection@\line@ \wdz@=\z@ \ht\z@=\z@ \dp\z@=\z@
	{\vDirection@(1,-1)\xydashl@ \xyatipfont\char\DirectionChar}%
	{\vDirection@(1,+1)\xydashl@ \xybtipfont\char\DirectionChar}%
}

\xydef@\ar@form{
	\ifx \space@\next \expandafter\DN@\space{\xyFN@\ar@form}%
	\else\ifx ^\next \DN@ ^{\xyFN@\ar@style}\edef\arvariant@@{\string^}%
	\else\ifx _\next \DN@ _{\xyFN@\ar@style}\edef\arvariant@@{\string_}%
	\else\ifx 0\next \DN@ 0{\xyFN@\ar@style}\def\arvariant@@{0}%
	\else\ifx 1\next \DN@ 1{\xyFN@\ar@style}\def\arvariant@@{1}%
	\else\ifx 2\next \DN@ 2{\xyFN@\ar@style}\def\arvariant@@{2}%
	\else\ifx 3\next \DN@ 3{\xyFN@\ar@style}\def\arvariant@@{3}%
	\else\ifx 4\next \DN@ 4{\xyFN@\ar@style}\def\arvariant@@{4}%
	\else\ifx \bgroup\next \let\next@=\ar@style
	\else\ifx [\next \DN@[##1]{\ar@modifiers{[##1]}}
	\else\ifx *\next \DN@ *{\ar@modifiers}%
	\else\addLT@\ifx\next \let\next@=\ar@slide
	\else\ifx /\next \let\next@=\ar@curveslash
	\else\ifx (\next \let\next@=\ar@curveinout 
	\else\addRQ@\ifx\next \addRQ@\DN@{\ar@curve@}%
	\else\addLQ@\ifx\next \addLQ@\DN@{\xyFN@\ar@curve}%
	\else\addDASH@\ifx\next \addDASH@\DN@{\defarstem@-\xyFN@\ar@}%
	\else\addEQ@\ifx\next \addEQ@\DN@{\def\arvariant@@{2}\defarstem@-\xyFN@\ar@}%
	\else\addDOT@\ifx\next \addDOT@\DN@{\defarstem@.\xyFN@\ar@}%
	\else\ifx :\next \DN@:{\def\arvariant@@{2}\defarstem@.\xyFN@\ar@}%
	\else\ifx ~\next \DN@~{\defarstem@~\xyFN@\ar@}%
	\else\ifx !\next \DN@!{\dasharstem@\xyFN@\ar@}%
	\else\ifx ?\next \DN@?{\ar@upsidedown\xyFN@\ar@}%
	\else \let\next@=\ar@error
	\fi\fi\fi\fi\fi\fi\fi\fi\fi\fi\fi\fi\fi\fi\fi\fi\fi\fi\fi\fi\fi\fi\fi \next@}

\makeatother

\newcommand{\Tone}
{\xymatrix@R=0.3cm@C=0cm{
	1\ot (A \prt B) \ar[rr]^{\distl} \ar[ddr]_{\lambda} && (1\ot A)\prt B \ar[ddl]^{\lambda \partimes \id} \\
	& {\lower -5pt\hbox{(T1)}} & \\
	& A\prt B & \\
}\label{eq:T1}
}

\newcommand{\Ttwo}
{\xymatrix@R=0.3cm@C=0cm{
		(A \prt B) \ot 1 \ar[rr]^{\distr} \ar[ddr]_{\ro} && A \prt {(B\otimes 1)} \ar[ddl]^{\id \prt \ro} \\
		& {\lower -5pt\hbox{(T2)}} & \\
		& A\prt B & \\
	}\label{eq:T2} 
}

\newcommand{\Tthree}
{\xymatrix@R=0.3cm@C=0cm{
		(K \prt A) \ot B \ar[rr]^{\distr} \ar[ddr]_{\lp \,\ot\, \id} && K \prt {(A \ot B)} \ar[ddl]^{\lp} \\
		& {\lower -5pt\hbox{(T3)}} & \\
		& A\ot B & \\
	}\label{eq:T3} 
}

\newcommand{\Tfour}
{\xymatrix@R=0.3cm@C=0cm{
		A \ot (B \prt K) \ar[rr]^{\distl} \ar[ddr]_{\id \,\ot\, \rp} && {(A\ot B)} \prt K \ar[ddl]^{\rp} \\
		& {\lower -5pt\hbox{(T4)}} & \\
		& A\ot B & \\
	}\label{eq:T4} 
}

\deftwocell[black]{ot : 2 -> 1}
\deftwocell[black,circle]{oteta : 0 -> 1}
\deftwocell[red]{par : 2 -> 1}
\deftwocell[red,circle]{pareta : 0 -> 1}
\deftwocell[crossing]{tau : 2 -> 2}
\deftwocell[gray]{alpha : 3 -> 1}
\deftwocell[gray,lefthalfcircle]{lambda : 1 -> 1}
\deftwocell[gray,righthalfcircle]{rho : 1 -> 1}
\deftwocell[gray,rectangle]{beta : 2 -> 1}
\deftwocell[gray,rectangle]{gamma : 3 -> 1}
\deftwocell[white]{aleph : 4 -> 1}
\deftwocell[white]{beth : 2 -> 1}
\deftwocell[white,rectangle]{gimmel : 2 -> 1}
\deftwocell[white,rectangle]{daleth : 3 -> 1}
\deftwocell[black]{id : 1 -> 1}

\deftwocell[dots]{dots : 2 -> 2}
\deftwocell[gray,rectangle]{f : 2 -> 2}

\tikzset{anchorbase/.style={baseline={([yshift=-0.5ex]current bounding box.center)}},
	tinynodes/.style={font=\tiny,text height=0.75ex,text depth=0.15ex},
	smallnodes/.style={font=\scriptsize,text height=0.75ex,text depth=0.15ex},
	>={Latex[length=1mm, width=1.5mm]}
}

\newcommand{\aosurf}{
	\adjustbox{margin=0.05cm}{
	\begin{tikzpicture} [anchorbase,xscale=-.2,yscale=.2,tinynodes]
		\draw [thick] (1,2) to [out=10,in=160] (7,1.5);
		\draw [thick] (2,1) to [out=10,in=150] (5,.5);
		\draw [thick] (3,0) to [out=0,in=210] (5,.5);
		\draw [thick] (5,.5) to [out=0,in=260] (7,1.5);
		\draw [thick] (7,1.5) to (8,1.5);
		\draw[thick]  (1,2) to (1,2+4);
		\draw[thick]  (2,1) to (2,1+4);
		\draw[thick]  (3,0) to (3,0+4);
		\draw[thick]  (8,1.5) to (8,1.5+4);
		\draw [thick] (1,2+4) to [out=10,in=160] (4,2+4);
		\draw [thick] (2,1+4) to [out=0,in=230] (4,2+4);
		\draw [thick] (3,0+4) to [out=0,in=230] (7,1.5+4);
		\draw [thick] (4,2+4) to[out=0,in=150] (7,1.5+4);
		\draw [thick] (7,1.5+4) to (8,1.5+4);
		\draw [line width=.6mm,] (5,.5) to [out=80,in=210] (5.75,0+2.85);
		\draw [line width=.6mm,] (7,1.5) to [out=90,in=-10] (5.75,0+2.85);
		\draw [line width=.6mm,] (5.75,0+2.85) to [out=30,in=270] (7,1.5+4);
		\draw [line width=.6mm,] (5.75,0+2.85) to [out=160,in=270] (4,2+4);
		\filldraw[draw=Purple, fill=Purple, ultra thick] (5.75,0+2.85) circle (0.3cm);
	\end{tikzpicture}
}
}

\newcommand{\apsurf}{
\adjustbox{margin=0.05cm}{
	\begin{tikzpicture} [anchorbase,xscale=-.2,yscale=.2,tinynodes]
		\draw [thick] (1,2) to [out=10,in=160] (7,1.5);
		\draw [thick] (2,1) to [out=10,in=150] (5,.5);
		\draw [thick] (3,0) to [out=0,in=210] (5,.5);
		\draw [thick] (5,.5) to [out=0,in=260] (7,1.5);
		\draw [thick] (7,1.5) to (8,1.5);
		\draw[thick]  (1,2) to (1,2+4);
		\draw[thick]  (2,1) to (2,1+4);
		\draw[thick]  (3,0) to (3,0+4);
		\draw[thick]  (8,1.5) to (8,1.5+4);
		\draw [thick] (1,2+4) to [out=10,in=160] (4,2+4);
		\draw [thick] (2,1+4) to [out=0,in=230] (4,2+4);
		\draw [thick] (3,0+4) to [out=0,in=230] (7,1.5+4);
		\draw [thick] (4,2+4) to[out=0,in=150] (7,1.5+4);
		\draw [thick] (7,1.5+4) to (8,1.5+4);
		\draw [line width=.6mm,color=red] (5,.5) to [out=80,in=210] (5.75,0+2.85);
		\draw [line width=.6mm,color=red] (7,1.5) to [out=90,in=-10] (5.75,0+2.85);
		\draw [line width=.6mm,color=red] (5.75,0+2.85) to [out=30,in=270] (7,1.5+4);
		\draw [line width=.6mm,color=red] (5.75,0+2.85) to [out=160,in=270] (4,2+4);
		\filldraw[draw=Purple, fill=Purple, ultra thick] (5.75,0+2.85) circle (0.3cm);
	\end{tikzpicture}
}
}

\newcommand{\distlsurf}{
	\adjustbox{margin=0.05cm}{
		\begin{tikzpicture} [anchorbase,xscale=-.2,yscale=.2,tinynodes]
			\draw [thick] (1,2) to [out=10,in=160] (7,1.5);
			\draw [thick] (2,1) to [out=10,in=150] (5,.5);
			\draw [thick] (3,0) to [out=0,in=210] (5,.5);
			\draw [thick] (5,.5) to [out=0,in=260] (7,1.5);
			\draw [thick] (7,1.5) to (8,1.5);
			\draw[thick]  (1,2) to (1,2+4);
			\draw[thick]  (2,1) to (2,1+4);
			\draw[thick]  (3,0) to (3,0+4);
			\draw[thick]  (8,1.5) to (8,1.5+4);
			\draw [thick] (1,2+4) to [out=10,in=160] (4,2+4);
			\draw [thick] (2,1+4) to [out=0,in=230] (4,2+4);
			\draw [thick] (3,0+4) to [out=0,in=230] (7,1.5+4);
			\draw [thick] (4,2+4) to[out=0,in=150] (7,1.5+4);
			\draw [thick] (7,1.5+4) to (8,1.5+4);
			\draw [line width=.6mm,color=black] (5,.5) to [out=80,in=210] (5.75,0+2.85);
			\draw [line width=.6mm,color=red] (7,1.5) to [out=90,in=-10] (5.75,0+2.85);
			\draw [line width=.6mm,color=black] (5.75,0+2.85) to [out=30,in=270] (7,1.5+4);
			\draw [line width=.6mm,color=red] (5.75,0+2.85) to [out=160,in=270] (4,2+4);
			\filldraw[draw=Purple, fill=Purple, ultra thick] (5.75,0+2.85) circle (0.3cm);
		\end{tikzpicture}
	}
}

\newcommand{\distrsurf}{
\adjustbox{margin=0.05cm}{
	\begin{tikzpicture} [anchorbase,xscale=-.2,yscale=.2,tinynodes]
		\draw [thick] (1,2) to [out=10,in=160] (7,1.5);
		\draw [thick] (2,1) to [out=10,in=150] (5,.5);
		\draw [thick] (3,0) to [out=0,in=210] (5,.5);
		\draw [thick] (5,.5) to [out=0,in=260] (7,1.5);
		\draw [thick] (7,1.5) to (8,1.5);
		\draw[thick]  (1,2) to (1,2+4);
		\draw[thick]  (2,1) to (2,1+4);
		\draw[thick]  (3,0) to (3,0+4);
		\draw[thick]  (8,1.5) to (8,1.5+4);
		\draw [thick] (1,2+4) to [out=10,in=160] (4,2+4);
		\draw [thick] (2,1+4) to [out=0,in=230] (4,2+4);
		\draw [thick] (3,0+4) to [out=0,in=230] (7,1.5+4);
		\draw [thick] (4,2+4) to[out=0,in=150] (7,1.5+4);
		\draw [thick] (7,1.5+4) to (8,1.5+4);
		\draw [line width=.6mm,color=red] (5,.5) to [out=80,in=210] (5.75,0+2.85);
		\draw [line width=.6mm,color=black] (7,1.5) to [out=90,in=-10] (5.75,0+2.85);
		\draw [line width=.6mm,color=red] (5.75,0+2.85) to [out=30,in=270] (7,1.5+4);
		\draw [line width=.6mm,color=black] (5.75,0+2.85) to [out=160,in=270] (4,2+4);
		\filldraw[draw=Purple, fill=Purple, ultra thick] (5.75,0+2.85) circle (0.3cm);
	\end{tikzpicture}
}
}

\newcommand{\gammalsurf}{
	\adjustbox{margin=0.05cm}{
		\begin{tikzpicture} [anchorbase,xscale=-.2,yscale=.2,tinynodes]
			\draw [thick] (1,2) to [out=10,in=160] (4,1.5);
			\draw [thick] (4,1.5) to [out=10,in=160] (7,1.5);
			\draw [thick] (2,1) to [out=-10,in=-130] (4,1.5);
			\draw [thick] (3,0) to [out=0,in=240] (7,1.5);
			\draw [thick] (7,1.5) to (8,1.5);
			\draw [thick] (1,2+4) to [out=10,in=160] (7,1.5+4);
			\draw [thick] (3,0+4) to [out=0,in=230] (7,1.5+4);
			\draw [thick] (7,1.5+4) to (8,1.5+4);
			\draw [line width=.6mm,color=red] (4,1.5) to [out=80,in=210] (5.75,0+3.3);
			\draw [line width=.6mm,color=black] (7,1.5) to [out=90,in=-10] (5.75,0+3.3);
			\draw [line width=.6mm,color=red] (5.75,0+3.3) to [out=30,in=270] (7,1.5+4);
			\draw [line width=.6mm,color=black] (5.75,0+3.3) to [out=170,in=90] (2,1);
			\filldraw[draw=Purple, fill=Purple, ultra thick] (5.75,0+3.3) circle (0.3cm);
			\draw[thick]  (1,2) to (1,2+4);
			\draw[thick]  (3,0) to (3,0+4);
			\draw[thick]  (8,1.5) to (8,1.5+4);
		\end{tikzpicture}
	}
}

\newcommand{\gammarsurf}{
	\adjustbox{margin=0.05cm}{
		\begin{tikzpicture} [anchorbase,xscale=-.2,yscale=.2,tinynodes]
			\draw [thick] (1,2) to [out=10,in=160] (4,1.5);
			\draw [thick] (4,1.5) to [out=10,in=160] (7,1.5);
			\draw [thick] (2,1) to [out=-10,in=-130] (4,1.5);
			\draw [thick] (3,0) to [out=0,in=240] (7,1.5);
			\draw [thick] (7,1.5) to (8,1.5);
			\draw [thick] (1,2+4) to [out=10,in=160] (7,1.5+4);
			\draw [thick] (3,0+4) to [out=0,in=230] (7,1.5+4);
			\draw [thick] (7,1.5+4) to (8,1.5+4);
			\draw [line width=.6mm,color=black] (4,1.5) to [out=80,in=210] (5.75,0+3.3);
			\draw [line width=.6mm,color=red] (7,1.5) to [out=90,in=-10] (5.75,0+3.3);
			\draw [line width=.6mm,color=black] (5.75,0+3.3) to [out=30,in=270] (7,1.5+4);
			\draw [line width=.6mm,color=red] (5.75,0+3.3) to [out=170,in=90] (2,1);
			\filldraw[draw=Purple, fill=Purple, ultra thick] (5.75,0+3.3) circle (0.3cm);
			\draw[thick]  (1,2) to (1,2+4);
			\draw[thick]  (3,0) to (3,0+4);
			\draw[thick]  (8,1.5) to (8,1.5+4);
		\end{tikzpicture}
	}
}

\newcommand{\losurf}{
\adjustbox{margin=0.05cm}{
	\begin{tikzpicture} [anchorbase,xscale=-.2,yscale=.2,tinynodes]
		\draw [thick] (1,2) to [out=10,in=160] (6,1.5);
		\draw [thick] (4,.5) to [out=0,in=260] (6,1.5);
		\draw [thick] (6,1.5) to (8,1.65);
		\draw[thick]  (1,2) to (1,2+3);
		\draw[thick]  (8,1.65) to (8,1.5+4);
		\draw [thick] (1,1+4) to (8,1.5+4);
		\draw [line width=.6mm,color=black] (4,.5) to [out=90,in=210] (4.75,0+3.5);
		\draw [line width=.6mm,color=black] (6,1.5) to [out=90,in=-10] (4.75,0+3.5);
		\filldraw[draw=Purple, fill=Purple, ultra thick] (4.75,0+3.5) circle (0.3cm);
	\end{tikzpicture}
}
}

\newcommand{\rosurf}{
	\adjustbox{margin=0.05cm}{
	\begin{tikzpicture} [anchorbase,xscale=-.2,yscale=.2,tinynodes]
		\draw [thick] (1,0.25) to [out=0,in=230] (6,1.5);
		\draw [thick] (3,2) to [out=30,in=110] (6,1.5);
		\draw [thick] (6,1.5) to (8,1.65);
		\draw[thick]  (1,0.25) to (1,2+3);
		\draw[thick]  (8,1.65) to (8,1.5+4);
		\draw [thick] (1,1+4) to (8,1.5+4);
		\draw [line width=.6mm,color=black] (3,2) to [out=90,in=180] (4.75,0+3.75);
		\draw [line width=.6mm,color=black] (6,1.5) to [out=90,in=-10] (4.75,0+3.75);
		\filldraw[draw=Purple, fill=Purple, ultra thick] (4.75,0+3.75) circle (0.3cm);
	\end{tikzpicture}
}
}

\newcommand{\lpsurf}{
\adjustbox{margin=0.05cm}{
	\begin{tikzpicture} [anchorbase,xscale=-.2,yscale=.2,tinynodes]
		\draw [thick] (1,2) to [out=10,in=160] (6,1.5);
		\draw [thick] (4,.5) to [out=0,in=260] (6,1.5);
		\draw [thick] (6,1.5) to (8,1.65);
		\draw[thick]  (1,2) to (1,2+3);
		\draw[thick]  (8,1.65) to (8,1.5+4);
		\draw [thick] (1,1+4) to (8,1.5+4);
		\draw [line width=.6mm,color=red] (4,.5) to [out=90,in=210] (4.75,0+3.5);
		\draw [line width=.6mm,color=red] (6,1.5) to [out=90,in=-10] (4.75,0+3.5);
		\filldraw[draw=Purple, fill=Purple, ultra thick] (4.75,0+3.5) circle (0.3cm);
	\end{tikzpicture}
}
}

\newcommand{\rpsurf}{
\adjustbox{margin=0.05cm}{
	\begin{tikzpicture} [anchorbase,xscale=-.2,yscale=.2,tinynodes]
		\draw [thick] (1,0.25) to [out=0,in=230] (6,1.5);
		\draw [thick] (3,2) to [out=30,in=110] (6,1.5);
		\draw [thick] (6,1.5) to (8,1.65);
		\draw[thick]  (1,0.25) to (1,2+3);
		\draw[thick]  (8,1.65) to (8,1.5+4);
		\draw [thick] (1,1+4) to (8,1.5+4);
		\draw [line width=.6mm,color=red] (3,2) to [out=90,in=180] (4.75,0+3.75);
		\draw [line width=.6mm,color=red] (6,1.5) to [out=90,in=-10] (4.75,0+3.75);
		\filldraw[draw=Purple, fill=Purple, ultra thick] (4.75,0+3.75) circle (0.3cm);
	\end{tikzpicture}
}
}

\newcommand{\whisker}{
	\scalebox{1.2}{
	\begin{tikzpicture} [anchorbase,xscale=-.2,yscale=.2,smallnodes]

		\draw [thick,lightgray] (17.5,3.95) to (17.5,-0.55);
		\draw [thick,lightgray] (17.5,-0.5) to (12.45,-1.9);
		
		\draw [thick,lightgray] (14.53,5.5) to (9.65,3.95);
		\draw [thick,lightgray] (14.5,5.5) to (14.5,0.8);
		\draw [thick,lightgray] (14.5,0.8) to (9.7,-0.6);
		\draw[thick,lightgray]  (9.7,-0.65) to (9.7,4);
		
		\draw[thick,lightgray]  (9.7,-0.6) to (12.5,-1.9);
		\draw[thick,lightgray]  (9.7,4) to (12.5,2.6);
		
		\draw[thick,lightgray]  (17.5,3.9) to (14.5,5.5);
		\draw [thick,lightgray] (17.5,-0.5) to (14.5,0.8);
		\node[lightgray] at (16, 2){$v_1$};
	
\draw [thick,lightgray] (2.5,11.5) to (-2.5,10.1);
\draw [thick,lightgray] (2.5,11.56) to (2.5,7);
\draw [thick,lightgray] (2.5,7) to (-2.5,5.6);
\draw [thick,lightgray] (-2.5,5.55) to (-2.5,10.15);

\draw [thick,lightgray] (5.5,9.95) to (5.5,5.15);
\draw [thick,lightgray] (5.5,5.2) to (0.7,3.8);

\draw [thick,lightgray] (0.7,3.8) to (-2.5,5.6);
\draw [thick,lightgray] (0.7,8.4) to (-2.5,10.1);
\draw [thick,lightgray] (5.5,9.9) to (2.5,11.5);
\draw [thick,lightgray] (5.5,5.2) to (2.5,7);

\node[lightgray] at (4.1,8){$w_1$};

	\draw [thick,lightgray] (12.5,-1.9) to (12.5,2.6);
	\draw [thick,lightgray] (17.5,3.9) to (12.5,2.6);
		
	\draw [thick,lightgray] (0.7,3.75) to (0.7,8.4);
	\draw [thick,lightgray] (5.5,9.9) to (0.7,8.4);
	
		\node[lightgray] at (-0.9, -2.6){$w_2$};
\draw [thick,lightgray] (10.5,9) to (25.5,1.5);
\draw [thick,lightgray] (10.5,13.5) to (25.5,6);
\draw[thick,lightgray]  (10.5,8.95) to (10.5,13.55);
\draw[thick,lightgray]  (25.5,1.45) to (25.5,6);
\draw [thick,lightgray] (29,7) to (14,14.5);
\draw [thick,lightgray] (14,14.5) to (10.5,13.5);
\draw [thick,lightgray] (29,7) to (25.5,6);
\draw [thick,lightgray] (29,7.05) to (29,2.45);
\draw [thick,lightgray] (29,2.5) to (25.5,1.5);
\draw[thick,lightgray]  (14,10) to (14,14.5);
\draw [thick,lightgray] (10.5,9) to (14,10);
\draw [thick,lightgray] (14,10) to (29,2.5);
		
	\node[lightgray] at (27.6, 4.9){$v_2$};
	\draw [thick,lightgray] (-18,1.5) to (-3,-6);
	\draw [thick,lightgray] (-14.5,2.4) to (0.5,-5);
	\draw [thick,lightgray] (-18,6) to (-3,-1.5);
	\draw [thick,lightgray]  (-18,1.45) to (-18,6.05);
	\draw [thick,lightgray]  (-3,-6.05) to (-3,-1.45);
	\draw [thick,lightgray] (0.5,-0.5) to (-14.5,6.9);
	\draw [thick,lightgray] (-14.5,6.9) to (-18,6);
	\draw [thick,lightgray] (0.5,-0.5) to (-3,-1.5);
	\draw [thick,lightgray] (0.5,-0.45) to (0.5,-5.05);
	\draw [thick,lightgray] (0.5,-5) to (-3,-6);
	\draw [thick,lightgray]  (-14.5,2.35) to (-14.5,6.95);
	\draw [thick,lightgray] (-18,1.5) to (-14.5,2.4);
	
	\draw [thick] (4.1,2.2) to [out=30,in=170] (7.5,3);
	\draw [thick] (5.9,1.8) to [out=0,in=260] (7.5,3);
	\draw [thick] (7.5,3) to (9.3,3.5);
	\draw[thick]  (4.1,2.14) to (4.1,6.85);
	\draw[thick]  (9.3,3.45) to (9.3,8.26);
	\draw [thick] (4.1,6.8) to (9.3,8.2);
	\draw [line width=.6mm,color=red] (5.9,1.72) to [out=100,in=220] (6.25,5);
	\draw [line width=.6mm,color=black] (7.5,3) to [out=90,in=-10] (6.25,5);
	\filldraw[draw=Purple, fill=Purple, ultra thick] (6.25,5) circle (0.3cm) node[above=-1, Purple]{$x$};
	
	\draw[thick,->,gray] (0-22-4,0+5+2.4,0-10) -- (4-22-4,1+5+2.4,0-10) node[above=0.2]{$\star_1$};
	\draw[thick,->,gray] (0-26.05,7.38,0-10) -- (0-26.05,7.38+5,0-10) node[right=-1]{$\star_2$};
	\draw[thick,->,gray] (-28,5.5,-15) -- (0-22,0+5.1,0-10) node[right=3]{$\star_0$};

	\end{tikzpicture}
}
}

\newcommand{\Partords}{\scalebox{1}{
		\xymatrix{
			& { \twocell{ (oteta *0 1) *1 ot } }
			\\
			{ \twocell{ (oteta *0 1 *0 oteta) *1 (1 *0 ot) *1 ot} }
			\ar@3 [ur] ^-{f_1\eq r \,\star_{1}\, s_2(l)}
			\ar@3 [rr] _-{\vphantom{\rule{0pt}{1.6ex}} f_2\eq s_2(r) \,\star_{1}\, l} ^-{}="1"
			&& { \twocell{ (1 *0 oteta) *1 ot } }
		}
	}
}

\newcommand{\Partordb}{\scalebox{1}{
		\xymatrix{
			& {  \twocell{ (oteta *0 oteta) *1 ot } }
			\\
			{ \twocell{ (oteta *0 oteta *0 oteta) *1 (1 *0 ot) *1 ot} }
			\ar@3 [ur] ^-{g_1 \eq v_2 \,\star_1\, r \,\star_{1}\, s_2(l)}
			\ar@3 [rr] _-{\vphantom{\rule{0pt}{1.6ex}} g_2 \eq v_2 \,\star_1\, s_2(r) \,\star_{1}\, l} ^-{}="1"
			&& { \twocell{ (oteta *0 oteta) *1 ot } }\;.
		}
	}
}

\newcommand{\Critdistlronc}{\scalebox{0.92}{
		\xymatrix{
			& { \twocell{ ( 1 *0 oteta *0 1) *1 (1 *0 par) *1 ot } }
			\\
			{ \twocell{ ( 1 *0 oteta *0 1) *1 (ot *0 1) *1 par} }
			\ar@3 [ur] ^-{\distlc}
			\ar@3 [rr] _-{\vphantom{\rule{0pt}{1.2ex}} \roc} ^-{}="1"
			&& { \twocell{par} }\;,
		}
	}
}

\newcommand{\Critdistrrpnc}{\scalebox{0.92}{
		\xymatrix{
			& { \twocell{ ( 1 *0 pareta *0 1) *1 (1 *0 ot) *1 par } }
			\\
			{ \twocell{ ( 1 *0 pareta *0 1) *1 (par *0 1) *1 ot} }
			\ar@3 [ur] ^-{\distrc}
			\ar@3 [rr] _-{\vphantom{\rule{0pt}{1.6ex}} \rpc} ^-{}="1"
			&& { \twocell{ot} }\;.
		}
	}
}

\newcommand{\Critdistllo}{\scalebox{0.9}{
	\xymatrix{
		& { \twocell{ ( oteta *0 1 *0 1) *1 (1 *0 par) *1 ot} }
		\ar@[lightgray]@3[dr] ^-{\color{lightgray}\loc}
		\\
		{ \twocell{ ( oteta *0 1 *0 1) *1 (ot *0 1) *1 par } }
		\ar@3 [ur] ^-{\distlc}
		\ar@3 [rr] _-{ \loc \vphantom{\rule{0pt}{1.6ex}}} ^-{}="1"
		&& { \twocell{par} }
		\ar@0 "2,2"!<20pt,15pt>;"1"!<0pt,0pt> _*+{(1)}
	}\label{Crit1}
}
}

\newcommand{\Critdistrro}{\scalebox{0.9}{
		\xymatrix{
			& { \twocell{ ( 1 *0 1 *0 oteta) *1 (1 *0 ot) *1 par } }
			\ar@[lightgray]@3[dr] ^-{\color{lightgray}\roc}
			\\
			{ \twocell{ ( 1 *0 1 *0 oteta) *1 (par *0 1) *1 ot} }
			\ar@3 [ur] ^-{\distrc}
			\ar@3 [rr] _-{\roc  \vphantom{\rule{0pt}{1.2ex}}} ^-{}="1"
			&& { \twocell{par} }
			\ar@0 "2,2"!<20pt,15pt>;"1"!<0pt,0pt> _*+{(2)}
		}\label{Crit2}
	}
}

\newcommand{\Critdistrlp}{\scalebox{0.9}{
		\xymatrix{
			& { \twocell{ ( pareta *0 1 *0 1) *1 (1 *0 ot) *1 par } }
			\ar@[lightgray]@3[dr] ^-{\color{lightgray}\lpc}
			\\
			{ \twocell{ ( pareta *0 1 *0 1) *1 (par *0 1) *1 ot} }
			\ar@3 [ur] ^-{\distrc}
			\ar@3 [rr] _-{\lpc  \vphantom{\rule{0pt}{1.8ex}}} ^-{}="1"
			&& { \twocell{ot} }
			\ar@0 "2,2"!<20pt,15pt>;"1"!<0pt,0pt> _*+{(3)}
		}\label{Crit3}
	}
}

\newcommand{\Critdistlrp}{\scalebox{0.9}{
		\xymatrix{
			& { \twocell{ ( 1 *0 1 *0 pareta) *1 (1 *0 par) *1 ot } }
			\ar@[lightgray]@3[dr] ^-{\color{lightgray}\rpc}
			\\
			{ \twocell{ ( 1 *0 1 *0 pareta) *1 (ot *0 1) *1 par} }
			\ar@3 [ur] ^-{\distlc}
			\ar@3 [rr] _-{\rpc  \vphantom{\rule{0pt}{1.6ex}}} ^-{}="1"
			&& { \twocell{ot} }
			\ar@0 "2,2"!<20pt,15pt>;"1"!<0pt,0pt> _*+{(4)}
		}\label{Crit4}
	}
}

\newcommand{\Critdistlao}{\scalebox{0.9}{
	\xymatrix@R=10pt@C=15pt{
		&{ \twocell{ (1 *0 ot *0 1) *1 (ot *0 1) *1 par } }
		\ar@[lightgray]@3 [rr] 
		^-{\color{lightgray}\distlc} 
		_-{}="1" 
		&& { \twocell{ (1 *0 ot *0 1) *1 (1 *0 par) *1 ot } }
		\ar@[lightgray]@3 [dr] ^-{\color{lightgray}\distlc}
		\\
		{ \twocell{(ot *0 2) *1 (ot *0 1) *1 par } }
		\ar@3 [ur] ^-{\aoc}
		\ar@3 [drr] _-{\distlc  \vphantom{\rule{0pt}{1.7ex}}}
		&&&& { \twocell{ (1 *0 1 *0 par) *1 (1 *0 ot) *1 ot } }
		\\
		&& { \twocell{ (ot *0 par) *1 ot }  }
		\ar@[lightgray]@3 [urr] _-{\color{lightgray}\aoc}
		\ar@0 "3,2"!<25pt,20pt>;"1"!<0pt,0pt> _-*+{(5)}
	}\label{Crit5}
}
}

\newcommand{\Critaodistl}{\scalebox{0.9}{
		\xymatrix@R=10pt@C=15pt{
		&{ \twocell{ (1 *0 par *0 1) *1 (ot *0 1) *1 ot } }
		\ar@[lightgray]@3 [rr] 
		^-{\color{lightgray}\aoc} 
		_-{}="1" 
		&& { \twocell{ (1 *0 par *0 1) *1 (1 *0 ot) *1 ot } }
		\ar@[lightgray]@3 [dr] ^-{\color{lightgray}\distrc}
		\\
		{ \twocell{(ot *0 2) *1 (par *0 1) *1 ot } }
		\ar@3 [ur] ^-{\distlc}
		\ar@3 [drr] _-{\distrc}
		&&&& { \twocell{ (1 *0 1 *0 ot) *1 (1 *0 par) *1 ot } }
		\\
		&& { \twocell{ (ot *0 ot) *1 par }  }
		\ar@[lightgray]@3 [urr] _-{\color{lightgray}\distlc}
		\ar@0 "3,2"!<25pt,20pt>;"1"!<0pt,0pt> _-*+{(6)}\label{Crit6}
	}
	}
}

\newcommand{\Critdistlap}{\scalebox{0.9}{
		\xymatrix@R=10pt@C=15pt{
			&{ \twocell{(1 *0 par *0 1) *1 (ot *0 1) *1 par } }
			\ar@[lightgray]@3 [rr] 
			^-{\color{lightgray}\distlc} 
			_-{}="1" 
			&& { \twocell{ (1 *0 par *0 1) *1 (1 *0 par) *1 ot } }
			\ar@[lightgray]@3 [dr] ^-{\color{lightgray}\apc}
			\\
			{ \twocell{ (ot *0 2) *1 (par *0 1) *1 par } }
			\ar@3 [ur] ^-{\distlc}
			\ar@3 [drr] _-{\apc}
			&&&& { \twocell{ (2 *0 par) *1 (1 *0 par) *1 ot } }
			\\
			&& { \twocell{ (ot *0 par) *1 par } }
			\ar@[lightgray]@3 [urr] _-{\color{lightgray}\distlc}
			\ar@0 "3,2"!<25pt,20pt>;"1"!<0pt,0pt> _-*+{(7)}
		}\label{Crit7}
	}
}

\newcommand{\Critdistrao}{\scalebox{0.9}{
		\xymatrix@R=10pt@C=15pt{
			&{ \twocell{(1 *0 ot *0 1) *1 (par *0 1) *1 ot } }
			\ar@[lightgray]@3 [rr] 
			^-{\color{lightgray}\distrc} 
			_-{}="1" 
			&& { \twocell{ (1 *0 ot *0 1) *1 (1 *0 ot) *1 par } }
			\ar@[lightgray]@3 [dr] ^-{\color{lightgray}\aoc}
			\\
			{ \twocell{ (par *0 2) *1 (ot *0 1) *1 ot } }
			\ar@3 [ur] ^-{\distrc}
			\ar@3 [drr] _-{\aoc}
			&&&& { \twocell{ (2 *0 ot) *1 (1 *0 ot) *1 par } }
			\\
			&& { \twocell{ (par *0 ot) *1 ot } }
			\ar@[lightgray]@3 [urr] _-{\color{lightgray}\distrc}
			\ar@0 "3,2"!<25pt,20pt>;"1"!<0pt,0pt> _-*+{(8)}
		}\label{Crit8}
	}
}

\newcommand{\Critdistldistr}{\scalebox{0.9}{
		\xymatrix@R=10pt@C=15pt{
			&{ \twocell{(1 *0 ot *0 1) *1 (par *0 1) *1 par } }
			\ar@[lightgray]@3 [rr] 
			^-{\color{lightgray}\apc} 
			_-{}="1" 
			&& { \twocell{ (1 *0 ot *0 1) *1 (1 *0 par) *1 par } }
			\ar@[lightgray]@3 [dr] ^-{\color{lightgray}\distrc}
			\\
			{ \twocell{(par *0 2) *1 (ot *0 1) *1 par } }
			\ar@3 [ur] ^-{\distrc}
			\ar@3 [drr] _-{\distlc}
			&&&& {  \twocell{ (2 *0 par) *1 (1 *0 ot) *1 par }  }
			\\
			&& {  \twocell{ (par *0 par) *1 ot } }
			\ar@[lightgray]@3 [urr] _-{\color{lightgray}\distrc}
			\ar@0 "3,2"!<25pt,20pt>;"1"!<0pt,0pt> _-*+{(9)}
		}\label{Crit9}
	}
}

\newcommand{\Critapdistr}{\scalebox{0.9}{
		\xymatrix@R=10pt@C=15pt{
			&{ \twocell{(1 *0 par *0 1) *1 (par *0 1) *1 ot }}
			\ar@[lightgray]@3 [rr] 
			^-{\color{lightgray}\distrc} 
			_-{}="1" 
			&& { \twocell{ (1 *0 par *0 1) *1 (1 *0 ot) *1 par } }
			\ar@[lightgray]@3 [dr] ^-{\color{lightgray}\distrc}
			\\
			{ \twocell{(par *0 2) *1 (par *0 1) *1 ot } }
			\ar@3 [ur] ^-{\apc}
			\ar@3 [drr] _-{\distrc}
			&&&& { \twocell{ (2 *0 ot) *1 (1 *0 par) *1 par } }
			\\
			&& { \twocell{ (par *0 ot) *1 par } }
			\ar@[lightgray]@3 [urr] _-{\color{lightgray}\apc}
			\ar@0 "3,2"!<20pt,20pt>;"1"!<0pt,0pt> _-*+{(10)}
		}\label{Crit10}
	}
}

\newcommand{\Critaoao}{\scalebox{0.9}{
		\xymatrix@R=10pt@C=15pt{
			&{ \twocell{(1 *0 ot *0 1) *1 (ot *0 1) *1 ot }}
			\ar@[lightgray]@3 [rr] 
			^-{\color{lightgray}\aoc} 
			_-{}="1" 
			&& { \twocell{ (1 *0 ot *0 1) *1 (1 *0 ot) *1 ot } }
			\ar@[lightgray]@3 [dr] ^-{\color{lightgray}\aoc}
			\\
			{ \twocell{(ot *0 2) *1 (ot *0 1) *1 ot } }
			\ar@3 [ur] ^-{\aoc}
			\ar@3 [drr] _-{\aoc}
			&&&& { \twocell{ (2 *0 ot) *1 (1 *0 ot) *1 ot } }
			\\
			&& { \twocell{ (ot *0 ot) *1 ot } }
			\ar@[lightgray]@3 [urr] _-{\color{lightgray}\aoc}
			\ar@0 "3,2"!<20pt,20pt>;"1"!<0pt,0pt> _-*+{(15)}
		}\label{Crit15}
	}
}

\newcommand{\Critaoro}{\scalebox{0.9}{
		\xymatrix{
			& { \twocell{ ( 1 *0 oteta *0 1) *1 (1 *0 ot) *1 ot } }
			\ar@[lightgray]@3[dr] ^-{\color{lightgray}\loc}
			\\
			{ \twocell{ ( 1 *0 oteta *0 1) *1 (ot *0 1) *1 ot} }
			\ar@3 [ur] ^-{\aoc}
			\ar@3 [rr] _-{\vphantom{\rule{0pt}{1.6ex}} \roc} ^-{}="1"
			&& { \twocell{ot} }
			\ar@0 "2,2"!<22pt,15pt>;"1"!<0pt,0pt> _*+{(11)}
		}\label{Crit11}
	}
}

\newcommand{\Critaolo}{\scalebox{0.9}{
		\xymatrix{
			& { \twocell{ ( oteta *0 1 *0 1) *1 (1 *0 ot) *1 ot } }
			\ar@[lightgray]@3[dr] ^-{\color{lightgray}\loc}
			\\
			{ \twocell{ ( oteta *0 1 *0 1) *1 (ot *0 1) *1 ot} }
			\ar@3 [ur] ^-{\aoc}
			\ar@3 [rr] _-{\vphantom{\rule{0pt}{1.6ex}} \loc} ^-{}="1"
			&& { \twocell{ot} }
			\ar@0 "2,2"!<22pt,15pt>;"1"!<0pt,0pt> _*+{(12)}
		}\label{Crit12}
	}
}

\newcommand{\Critaororo}{\scalebox{0.9}{
		\xymatrix{
			& { \twocell{ ( 1 *0 1 *0 oteta) *1 (1 *0 ot) *1 ot } }
			\ar@[lightgray]@3[dr] ^-{\color{lightgray}\roc}
			\\
			{ \twocell{ ( 1 *0 1 *0 oteta) *1 (ot *0 1) *1 ot} }
			\ar@3 [ur] ^-{\aoc}
			\ar@3 [rr] _-{\vphantom{\rule{0pt}{1.6ex}} \roc} ^-{}="1"
			&& { \twocell{ot} }
			\ar@0 "2,2"!<22pt,15pt>;"1"!<0pt,0pt> _*+{(13)}
		}\label{Crit13}
	}
}

\newcommand{\Critloro}{\scalebox{0.9}{
	\xymatrix{
		{\twocell{(oteta *0 oteta) *1 ot}}
		\ar@3@/^5ex/ [rr] ^{\vphantom{\rule{0pt}{1.2ex}} \loc} _{}="1"
		\ar@3@/_5ex/ [rr] _{\vphantom{\rule{0pt}{1.2ex}} \roc} ^{}="2"
		&& {\twocell{oteta}}
		\ar@0 "1,1"!<-4pt,-7pt>;"1"!<20pt,15pt> _-*+{(14)}
	}\label{Crit14}
}
}

\newcommand{\Critapap}{\scalebox{0.9}{
		\xymatrix@R=10pt@C=15pt{
			&{ \twocell{(1 *0 par *0 1) *1 (par *0 1) *1 par }}
			\ar@[lightgray]@3 [rr] 
			^-{\color{lightgray}\apc} 
			_-{}="1" 
			&& { \twocell{ (1 *0 par *0 1) *1 (1 *0 par) *1 par } }
			\ar@[lightgray]@3 [dr] ^-{\color{lightgray}\apc}
			\\
			{ \twocell{(par *0 2) *1 (par *0 1) *1 par } }
			\ar@3 [ur] ^-{\apc}
			\ar@3 [drr] _-{\vphantom{\rule{0pt}{1.4ex}}\apc}
			&&&& { \twocell{ (2 *0 par) *1 (1 *0 par) *1 par } }
			\\
			&& { \twocell{ (par *0 par) *1 par } }
			\ar@[lightgray]@3 [urr] _-{\color{lightgray}\apc}
			\ar@0 "3,2"!<20pt,20pt>;"1"!<0pt,0pt> _-*+{(20)}
		}\label{Crit20}
	}
}

\newcommand{\Critaprp}{\scalebox{0.9}{
		\xymatrix{
			& { \twocell{ ( 1 *0 pareta *0 1) *1 (1 *0 par) *1 par } }
			\ar@[lightgray]@3[dr] ^-{\color{lightgray}\lpc}
			\\
			{ \twocell{ ( 1 *0 pareta *0 1) *1 (par *0 1) *1 par} }
			\ar@3 [ur] ^-{\apc}
			\ar@3 [rr] _-{\vphantom{\rule{0pt}{1.4ex}} \rpc} ^-{}="1"
			&& { \twocell{par} }
			\ar@0 "2,2"!<20pt,15pt>;"1"!<0pt,0pt> _*+{(16)}
		}\label{Crit16}
	}
}

\newcommand{\Critaplp}{\scalebox{0.9}{
		\xymatrix{
			& { \twocell{ ( pareta *0 1 *0 1) *1 (1 *0 par) *1 par } }
			\ar@[lightgray]@3[dr] ^-{\color{lightgray}\lpc}
			\\
			{ \twocell{ ( pareta *0 1 *0 1) *1 (par *0 1) *1 par} }
			\ar@3 [ur] ^-{\apc}
			\ar@3 [rr] _-{\vphantom{\rule{0pt}{1.5ex}} \lpc} ^-{}="1"
			&& { \twocell{par} }
			\ar@0 "2,2"!<20pt,15pt>;"1"!<0pt,0pt> _*+{(17)}
		}\label{Crit17}
	}
}

\newcommand{\Critaprprp}{\scalebox{0.9}{
		\xymatrix{
			& { \twocell{ ( 1 *0 1 *0 pareta) *1 (1 *0 par) *1 par } }
			\ar@[lightgray]@3[dr] ^-{\color{lightgray}\rpc}
			\\
			{ \twocell{ ( 1 *0 1 *0 pareta) *1 (par *0 1) *1 par} }
			\ar@3 [ur] ^-{\apc}
			\ar@3 [rr] _-{\vphantom{\rule{0pt}{1.6ex}} \rpc} ^-{}="1"
			&& { \twocell{par} }
			\ar@0 "2,2"!<20pt,15pt>;"1"!<0pt,0pt> _*+{(18)}
		}\label{Crit18}
	}
}

\newcommand{\Critlprp}{\scalebox{0.9}{
		\xymatrix{
			{\twocell{(pareta *0 pareta) *1 par}}
			\ar@3@/^5ex/ [rr] ^{\lpc} _{}="1"
			\ar@3@/_5ex/ [rr] _{\vphantom{\rule{0pt}{1.5ex}} \rpc} ^{}="2"
			&& {\twocell{pareta}}
			\ar@0 "1,1"!<-4pt,-7pt>;"1"!<20pt,15pt> _-*+{(19)}
		}\label{Crit19}
	}
}

\newcommand{\Critdistlroncr}{\scalebox{0.9}{
		\xymatrix{
			& { \twocell{ ( 1 *0 oteta *0 1) *1 (1 *0 par) *1 ot } }
			\ar@[lightgray]@3[dr] ^-{\color{lightgray}k^l}
			\\
			{ \twocell{ ( 1 *0 oteta *0 1) *1 (ot *0 1) *1 par} }
			\ar@3 [ur] ^-{\distlc}
			\ar@3 [rr] _-{\vphantom{\rule{0pt}{1.2ex}} \roc} ^-{}="1"
			&& { \twocell{par} }
			\ar@0 "2,2"!<20pt,15pt>;"1"!<0pt,0pt> _*+{(21)}
		}\label{Crit21}
	}
}

\newcommand{\Critdistrrpncr}{\scalebox{0.9}{
		\xymatrix{
			& { \twocell{ ( 1 *0 pareta *0 1) *1 (1 *0 ot) *1 par } }
			\ar@[lightgray]@3[dr] ^-{\color{lightgray}k^r}
			\\
			{ \twocell{ ( 1 *0 pareta *0 1) *1 (par *0 1) *1 ot} }
			\ar@3 [ur] ^-{\distrc}
			\ar@3 [rr] _-{\vphantom{\rule{0pt}{1.5ex}} \rpc} ^-{}="1"
			&& { \twocell{ot} }
			\ar@0 "2,2"!<20pt,15pt>;"1"!<0pt,0pt> _*+{(22)}
		}\label{Crit22}
	}
}

\newcommand{\Critglrp}{\scalebox{0.9}{
		\xymatrix{
			{ \twocell{ ( 1 *0 oteta *0 pareta) *1 (1 *0 par) *1 ot } }
			\ar@3[rr] ^-{k^l}
			\ar@3 [d] ^-{\hspace{2pt} \rpc}
			&& { \twocell{ (1 *0 pareta) *1 par} }
			\ar@[lightgray]@3 [d] ^-{\hspace{2pt} \color{lightgray}\rpc}
			\\
			{ \twocell{ (1 *0 oteta) *1 ot } }
			\ar@[lightgray]@3 [rr] _-{\vphantom{\rule{0pt}{1.5ex}} \color{lightgray}\roc} ^-{}="1"
			&& { \twocell{id *1 id} }
			\ar@0 "2,2"!<30pt,30pt>;"1"!<0pt,0pt> _*+{(23)}
		}\label{Crit23}
	}
}

\newcommand{\Critgllo}{\scalebox{0.9}{
	\xymatrix{
		{\twocell{ ( oteta *0 oteta *0 1) *1 (1 *0 par) *1 ot }}
		\ar@3@/^5ex/ [rr] ^{k^l} _{}="1"
		\ar@3@/_5ex/ [rr] _{\vphantom{\rule{0pt}{1.5ex}} \loc} ^{}="2"
		&& {\twocell{(oteta *0 1) *1 par }}
		\ar@0 "1,1"!<-1pt,-7pt>;"1"!<20pt,15pt> _-*+{(24)}
	}\label{Crit24}
	}
}

\newcommand{\Critgldistr}{\scalebox{0.9}{
		\xymatrix{
			{ \twocell{ ( 2 *0 oteta *0 1) *1 (par *0 par) *1 ot } }
			\ar@3[rr] ^-{k^l}
			\ar@3 [d] ^-{\hspace{2pt}\distrc}
			&& { \twocell{ (par *0 1) *1 par} }
			\ar@[lightgray]@3 [d] ^-{\hspace{2pt}\color{lightgray}\apc}
			\\
			{ \twocell{(2 *0 oteta *0 1) *1 (2 *0 par) *1 (1 *0 ot) *1 par } }
			\ar@[lightgray]@3 [rr] _-{\vphantom{\rule{0pt}{1.5ex}} \color{lightgray}k^l} ^-{}="1"
			&& { \twocell{ (1 *0 par) *1 par } }
			\ar@0 "2,2"!<30pt,30pt>;"1"!<0pt,0pt> _*+{(25)}
		}\label{Crit25}
	}
}

\newcommand{\Critglao}{\scalebox{0.9}{
		\xymatrix{
			{ \twocell{ ( 2 *0 oteta *0 1) *1 (ot *0 par) *1 ot } }
			\ar@3[rr] ^-{k^l}
			\ar@3 [d] ^-{\hspace{2pt} \aoc}
			&& { \twocell{ (ot *0 1) *1 par} }
			\ar@[lightgray]@3 [d] ^-{\hspace{2pt} \color{lightgray}\distlc}
			\\
			{ \twocell{(2 *0 oteta *0 1) *1 (2 *0 par) *1 (1 *0 ot) *1 ot } }
			\ar@[lightgray]@3 [rr] _-{\vphantom{\rule{0pt}{1.5ex}} \color{lightgray}k^l} ^-{}="1"
			&& { \twocell{ (1 *0 par) *1 ot } }
			\ar@0 "2,2"!<30pt,30pt>;"1"!<0pt,0pt> _*+{(26)}
		}\label{Crit26}
	}
}

\newcommand{\Critgldistl}{\scalebox{0.9}{
		\xymatrix@R=10pt@C=15pt{
			&{ \twocell{(1 *0 oteta *0 2) *1 (1 *0 par *0 1) *1 (1 *0 par) *1 ot }}
			\ar@[lightgray]@3 [rr] 
			^-{\color{lightgray}\apc} 
			_-{}="1" 
			&& { \twocell{ (1 *0 oteta *0 par) *1 (1 *0 par) *1 ot } }
			\ar@[lightgray]@3 [dr] ^-{\color{lightgray}k^l}
			\\
			{ \twocell{ (1 *0 oteta *0 2) *1 (1 *0 par *0 1) *1 (ot *0 1) *1 par } }
			\ar@3 [ur] ^-{\distlc}
			\ar@3 [drr] _-{\vphantom{\rule{0pt}{1.7ex}} k^l}
			&&&& { \twocell{ (1 *0 par) *1 par } }
			\\
			&& { \twocell{ (par *0 1) *1 par } }
			\ar@[lightgray]@3 [urr] _-{\color{lightgray}\apc}
			\ar@0 "3,2"!<20pt,20pt>;"1"!<0pt,0pt> _-*+{(27)}
		}\label{Crit27}
	}
}

\newcommand{\CritglaoP}{\scalebox{0.9}{
		\xymatrix@R=10pt@C=15pt{
			&{ \twocell{(1 *0 oteta *0 2) *1 (1 *0 par *0 1) *1 (1 *0 ot) *1 ot }}
			\ar@[lightgray]@3 [rr] 
			^-{\color{lightgray}\distrc} 
			_-{}="1" 
			&& { \twocell{ (1 *0 oteta *0 ot) *1 (1 *0 par) *1 ot } }
			\ar@[lightgray]@3 [dr] ^-{\color{lightgray}k^l}
			\\
			{ \twocell{ (1 *0 oteta *0 2) *1 (1 *0 par *0 1) *1 (ot *0 1) *1 ot } }
			\ar@3 [ur] ^-{\aoc}
			\ar@3 [drr] _-{\vphantom{\rule{0pt}{1.7ex}} k^l}
			&&&& { \twocell{ (1 *0 ot) *1 par } }
			\\
			&& { \twocell{ (par *0 1) *1 ot } }
			\ar@[lightgray]@3 [urr] _-{\color{lightgray}\distrc}
			\ar@0 "3,2"!<20pt,20pt>;"1"!<0pt,0pt> _-*+{(28)}
		}\label{Crit28}
	}
}

\newcommand{\Critgrro}{\scalebox{0.9}{
		\xymatrix{
			{ \twocell{ ( 1 *0 pareta *0 oteta) *1 (1 *0 ot) *1 par } }
			\ar@3[rr] ^-{k^r}
			\ar@3 [d] ^-{\hspace{2pt} \roc}
			&& { \twocell{ (1 *0 oteta) *1 ot} }
			\ar@[lightgray]@3 [d] ^-{\hspace{2pt} \color{lightgray}\roc}
			\\
			{ \twocell{ (1 *0 pareta) *1 par } }
			\ar@[lightgray]@3 [rr] _-{\vphantom{\rule{0pt}{1.5ex}} \color{lightgray}\rpc} ^-{}="1"
			&& { \twocell{id *1 id} }
			\ar@0 "2,2"!<27pt,22pt>;"1"!<0pt,0pt>_*+{(29)}
		}\label{Crit29}
	}
}

\newcommand{\Critgrlp}{\scalebox{0.9}{
		\xymatrix{
			{\twocell{ ( pareta *0 pareta *0 1) *1 (1 *0 ot) *1 par }}
			\ar@3@/^5ex/ [rr] ^{k^r} _{}="1"
			\ar@3@/_5ex/ [rr] _{\vphantom{\rule{0pt}{1.5ex}} \lpc} ^{}="2"
			&& {\twocell{(pareta *0 1) *1 ot }}
			\ar@0 "1,1"!<0pt,-7pt>;"1"!<20pt,15pt> _-*+{(30)}
		}\label{Crit30}
	}
}

\newcommand{\Critgrdistl}{\scalebox{0.9}{
		\xymatrix{
			{ \twocell{ ( 2 *0 pareta *0 1) *1 (ot *0 ot) *1 par } }
			\ar@3[rr] ^-{k^r}
			\ar@3 [d] ^-{\hspace{2pt} \distlc}
			&& { \twocell{ (ot *0 1) *1 ot} }
			\ar@[lightgray]@3 [d] ^-{\hspace{2pt} \color{lightgray}\aoc}
			\\
			{ \twocell{(2 *0 pareta *0 1) *1 (2 *0 ot) *1 (1 *0 par) *1 ot } }
			\ar@[lightgray]@3 [rr] _-{\vphantom{\rule{0pt}{1.5ex}} \color{lightgray}k^r} ^-{}="1"
			&& { \twocell{ (1 *0 ot) *1 ot } }
			\ar@0 "2,2"!<30pt,30pt>;"1"!<0pt,0pt>_*+{(31)}
		}\label{Crit31}
	}
}

\newcommand{\Critgrap}{\scalebox{0.9}{
		\xymatrix{
			{ \twocell{ ( 2 *0 pareta *0 1) *1 (par *0 ot) *1 par } }
			\ar@3[rr] ^-{k^r}
			\ar@3 [d] ^-{\hspace{2pt} \apc}
			&& { \twocell{ (par *0 1) *1 ot} }
			\ar@[lightgray]@3 [d] ^-{\hspace{2pt} \color{lightgray}\distrc}
			\\
			{ \twocell{(2 *0 pareta *0 1) *1 (2 *0 ot) *1 (1 *0 par) *1 par } }
			\ar@[lightgray]@3 [rr] _-{\vphantom{\rule{0pt}{1.5ex}} \color{lightgray}k^r} ^-{}="1"
			&& { \twocell{ (1 *0 ot) *1 par } }
			\ar@0 "2,2"!<30pt,30pt>;"1"!<0pt,0pt>_*+{(32)}
		}\label{Crit32}
	}
}

\newcommand{\Critgrdistr}{\scalebox{0.9}{
		\xymatrix@R=10pt@C=15pt{
			&{ \twocell{(1 *0 pareta *0 2) *1 (1 *0 ot *0 1) *1 (1 *0 ot) *1 par }}
			\ar@[lightgray]@3 [rr] 
			^-{\color{lightgray}\aoc} 
			_-{}="1" 
			&& { \twocell{ (1 *0 pareta *0 ot) *1 (1 *0 ot) *1 par } }
			\ar@[lightgray]@3 [dr] ^-{\color{lightgray}k^r}
			\\
			{ \twocell{ (1 *0 pareta *0 2) *1 (1 *0 ot *0 1) *1 (par *0 1) *1 ot } }
			\ar@3 [ur] ^-{\distrc}
			\ar@3 [drr] _-{\vphantom{\rule{0pt}{1.5ex}} k^r}
			&&&& { \twocell{ (1 *0 ot) *1 ot } }
			\\
			&& { \twocell{ (ot *0 1) *1 ot } }
			\ar@[lightgray]@3 [urr] _-{\color{lightgray}\aoc}
			\ar@0 "3,2"!<20pt,20pt>;"1"!<0pt,0pt> _-*+{(33)}
		}\label{Crit33}
	}
}

\newcommand{\CritgrapP}{\scalebox{0.9}{
		\xymatrix@R=10pt@C=15pt{
			&{ \twocell{(1 *0 pareta *0 2) *1 (1 *0 ot *0 1) *1 (1 *0 par) *1 par }}
			\ar@[lightgray]@3 [rr] 
			^-{\color{lightgray}\distlc} 
			_-{}="1" 
			&& { \twocell{ (1 *0 pareta *0 par) *1 (1 *0 ot) *1 par } }
			\ar@[lightgray]@3 [dr] ^-{\color{lightgray}k^r}
			\\
			{ \twocell{ (1 *0 pareta *0 2) *1 (1 *0 ot *0 1) *1 (par *0 1) *1 par } }
			\ar@3 [ur] ^-{\apc}
			\ar@3 [drr] _-{\vphantom{\rule{0pt}{1.5ex}} k^r}
			&&&& { \twocell{ (1 *0 par) *1 ot } }
			\\
			&& { \twocell{ (ot *0 1) *1 par } }
			\ar@[lightgray]@3 [urr] _-{\color{lightgray}\distlc}
			\ar@0 "3,2"!<20pt,20pt>;"1"!<0pt,0pt> _-*+{(34)}
		}\label{Crit34}
	}
}

\newcommand{\Critglgro}{\scalebox{0.9}{
		\xymatrix{
			{ \twocell{ (1 *0 oteta *0 pareta *0 1) *1 ( 2 *0 ot) *1 (1 *0 par) *1 ot } }
			\ar@3[rr] ^-{k^l}
			\ar@3 [d] ^-{\hspace{2pt} k^r}
			&& { \twocell{ (1 *0 pareta *0 1) *1 (1 *0 ot) *1 par} }
			\ar@[lightgray]@3 [d] ^-{\hspace{2pt} \color{lightgray}k^r}
			\\
			{ \twocell{(1 *0 oteta *0 1) *1 (1 *0 ot) *1 ot } }
			\ar@[lightgray]@3 [rr] _-{\vphantom{\rule{0pt}{1.5ex}} \color{lightgray}\loc} ^-{}="1"
			&& { \twocell{ot} }
			\ar@0 "2,2"!<30pt,30pt>;"1"!<0pt,0pt>_*+{(35)}
		}\label{Crit35}
	}
}

\newcommand{\Critglgrt}{\scalebox{0.9}{
	\xymatrix{
		{ \twocell{ (1 *0 pareta *0 oteta *0 1) *1 ( 2 *0 par) *1 (1 *0 ot) *1 par } }
		\ar@3[rr] ^-{k^l}
		\ar@3 [d] ^-{\hspace{2pt} k^r}
		&& { \twocell{ (1 *0 pareta *0 1) *1 (1 *0 par) *1 par} }
		\ar@[lightgray]@3 [d] ^-{\hspace{2pt} \color{lightgray}\lpc}
		\\
		{ \twocell{(1 *0 oteta *0 1) *1 (1 *0 par) *1 ot } }
		\ar@[lightgray]@3 [rr] _-{\vphantom{\rule{0pt}{1.5ex}} \color{lightgray}k^l} ^-{}="1"
		&& { \twocell{par} }
				\ar@0 "2,2"!<30pt,30pt>;"1"!<0pt,0pt>_*+{(36)}
		}\label{Crit36}
	}
}

\title[Higher rewriting for strong linearly distributive coherence]{Higher rewriting for\\ strong linearly distributive coherence}
\author{Max Demirdilek, Christian Reiher, Christoph Schweigert}

\begin{document}

\begin{abstract}
	Rewriting is a computational method to solve word problems: deciding whether two compound expressions, or words, are equivalent with respect to a set of identities. Higher rewriting generalizes rewriting to higher-dimensional words, or cells, in higher categories.
	
	Following Guiraud and Malbos, we use higher rewriting to generalize Mac Lane's coherence theorem. Our result concerns strong linearly distributive categories: categories equipped with two unital tensor products related by invertible distributors. We show that in a strong linearly distributive category every formal diagram commutes.
\end{abstract}

\maketitle
\vspace{-0.85cm}
\renewcommand{\baselinestretch}{1.08}\normalsize
\tableofcontents
\renewcommand{\baselinestretch}{1}\normalsize
\vspace{-1.05cm}

\section{Introduction}
\addtocontents{toc}{\protect\setcounter{tocdepth}{1}}
\subsection{Rewriting} 
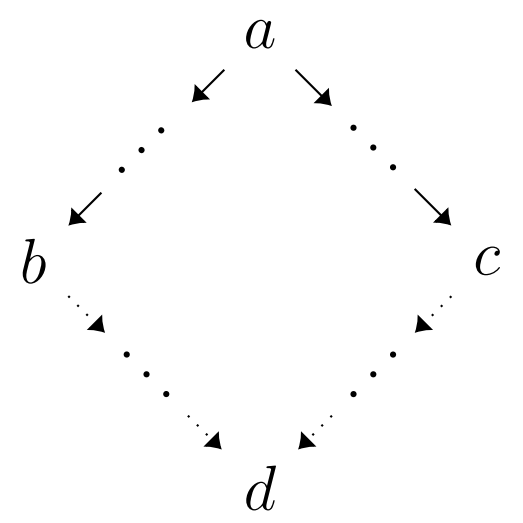
\begin{figure}[H]
	\begin{tikzpicture}[xscale=1,yscale=3]
		\node   (w) at (0,0) {$a$};
		\node[rotate=45] (xld) at (-0.75,-0.25) {$\cdots$};
		\node[rotate=-45]  (yld) at (0.75,-0.25) {$\ldots$};
		\node   (xm) at (-1.5,-0.5) {$b$};  
		\node   (yn) at (1.5,-0.5) {$c$};
		\node[rotate=-45]   (uld) at (-0.75,-0.75) {$\ldots$};
		\node[rotate=45]   (vld) at (0.75,-0.75) {$\ldots$};
		\node   (z) at (0,-1) {$d$};
		\draw[->] (w) to (xld);
		\draw[->] (w) to (yld);
		\draw[->] (xld) to (xm);
		\draw[->] (yld) to (yn);
		\draw[->,dotted] (xm) to (uld);
		\draw[->,dotted] (yn) to (vld);
		\draw[->,dotted] (vld) to (z);
		\draw[->,dotted] (uld) to (z);
	\end{tikzpicture}
	\caption{Confluence: any two 
		finite rewrite sequences $a \xrightarrow{\ast}b$ and $a \xrightarrow{\ast}c$ starting from the same word $a$ can be completed to a common word $d$.}
	\label{fig:diamond}
\end{figure}

Rewriting \cite{New42,BaNi98} has applications in computer science \cite{Hs85,DeJ90,RoHa05}, logic \cite{Me12}, computational algebra \cite{KnB70, MaRe98}, and linguistics \cite{Fr98}. It is a method to solve word problems: deciding whether two compound expressions are equivalent with respect to a set of identities. This is done by orienting these identities to obtain a binary relation $\to$ on the set of words. Such a \emph{rewrite system} is called \emph{terminating} if there does not exist an infinite \emph{rewrite sequence} $a_1\to a_2 \to \ldots $\,. It is called \emph{confluent} if any two finite rewrite sequences starting from the same word can be completed to rewrite sequences ending in the same word; see Figure \ref{fig:diamond}.

Terminating and confluent rewrite systems have unique normal forms. Their word problem can thus be solved by comparing normal forms, e.g. \cite[Thm.~2.1.9]{BaNi98}.
\subsection{Higher rewriting}\label{sec:Higher-rewriting}
Rewrite systems freely generate categories: words become objects, and rewrite sequences become non-identity morphisms. Higher rewriting \cite{ABGMMM25} is rewriting for higher-dimensional words, or cells, in higher categories, meaning words that can be composed in multiple ways. Higher rewriting studies rewrites between words, rewrites between rewrites, rewrites between rewrites between rewrites, et cetera. These data form a higher rewrite system called a \emph{computad}; see Definition \ref{def:computad}. 

\subsection{Linearly distributive coherence}
The purpose of this note is to illustrate, following \cite{GuMa,GuMa12}, how higher rewriting can be used to prove categorical coherence theorems: we apply higher rewriting to the coherence problem of linearly distributive categories. 

As recalled in Definition \ref{def:LD-cat}, a \emph{linearly distributive category} \cite{CS97} is a category equipped with two monoidal structures $(\otimes,1)$ and $(\partimes,K)$, with associators $\ao,\ap$ and unitors $\lo,\ro,\lp,\rp$. It is further endowed with two (not necessarily invertible) natural transformations
\begin{equation*}\begin{array}{lll}
		\distl_{A,B,C}:\: A\ot(B \prt C) & \longrightarrow &{(A \ot B)} \prt {C}\,,\\[0.4em]
		\distr_{A,B,C}:\: {(A \prt B)} \ot C & \longrightarrow &A \prt{(B \ot C)}\,,
	\end{array} \label{eq:distributors}
\end{equation*} 
 called \emph{distributors}. These are subject to numerous coherence conditions. Examples of linearly distributive categories include representation categories of vertex operator algebras \cite{ALSW25}, categories of topological vector spaces \cite{Ba00}, and derived categories of sheaves \cite{BoD13}. 

The linearly distributive \emph{coherence problem} asks: are any two parallel coherence morphisms (built from associators, unitors, and distributors, using composition and tensor products) equivalent with respect to the coherence conditions? It is known \cite{DP07, DRS26} that the answer is affirmative in a linearly distributive category \emph{without units}. If one includes the units, however, not every formal diagram commutes \cite{BCST96}; see also Examples \ref{ex:DRS counterexample} and \ref{ex:counterexample full coherence}. 

We use higher rewriting to prove that invertibility implies coherence with units:

\begin{maintheorem*}[Theorem \ref{thm:main theorem}]
	Linearly distributive categories with \emph{invertible} distributors are coherent, meaning all their formal diagrams commute.
\end{maintheorem*}

\subsection{Proof outline}
Definition \ref{def:3-computad LD} constructs a $3$-computad $\ms{LD}$ by orienting the coherence morphisms of \emph{strong} linearly distributive categories, meaning linearly distributive categories with \emph{invertible} distributors. The computad $\ms{LD}$ has a unique $1$-cell representing the underlying category, while its $2$-cells represent the monoidal products and units, and its $3$-cells represent the coherence morphisms. 

Subsection \ref{sec: conf} shows that completing $\ms{LD}$ to the $3$-computad $\highoverline{\ms{LD}}$ of Definition \ref{def: Knuth-Bendix completion LD} makes it \emph{confluent}, in the sense of Definition \ref{def: branch and confl}. The proof proceeds by computing the \emph{critical branchings} of $\highoverline{\ms{LD}}$; this notion is recalled in Definition \ref{def: critical branchings}. Subsection \ref{subs:termination} establishes that $\highoverline{\ms{LD}}$ is also \emph{terminating}; this notion is recalled in Definition \ref{def:termination}. 

By Proposition \ref{prop: homotop basis free track}, a result from the literature, the categorical structures presented by a confluent, terminating $3$-computad and by its critical branchings are coherent. To prove Theorem \ref{thm:main theorem}, we show that, for $\highoverline{\ms{LD}}$, these structures are strong linearly distributive categories.

\subsection{Alternative proof}
An alternative proof of Theorem \ref{thm:main theorem} uses that strong linearly distributive categories are suitably equivalent \cite{CS97,De22} to \emph{shift monoidal categories}: linearly distributive categories whose second monoidal product is ${X} \partimes {Y} ={X} \otimes {{K^{-1}} \otimes {Y}}$, where $K$ is an $\otimes$-invertible object with ${\otimes}$-inverse $K^{-1}$. Mac Lane's coherence theorem for monoidal categories can then be used to prove Theorem \ref{thm:main theorem}.

\addtocontents{toc}{\protect\setcounter{tocdepth}{2}}
\section{Strong linearly distributive categories}

\begin{definition}\label{def:LD-cat}
	A \emph{linearly distributive category} is a category $\cC$ equipped with two monoidal structures $(\otimes,1)$ and $(\partimes,K)$, and two (not necessarily invertible) natural transformations 
	\begin{align}\label{eq: distributorl}
		\distl\colon\; {\otimes} \circ (\idC \tim {\partimes}) &\;\slongrightarrow\; {\partimes} \circ ({\otimes} \tim \idC)\,,\\
		\distr\colon\; {\otimes} \circ ({\partimes} \tim \idC) &\;\slongrightarrow\; {\partimes} \circ ({\idC} \tim {\otimes})\,,
	\end{align}
	called \emph{distributors}, subject to six pentagon and four triangle axioms, where $A,B,C,D\in \cC$:
	
	\begingroup
	\renewcommand{\theequation}{\textup{P1}}
	\refstepcounter{equation}\label{eq:P1}
	\begin{equation*}
		\begin{tikzcd}[column sep=-1cm, row sep=0.28cm]
			A\ot (B\ot (C\partimes D)) && (A\ot B)\ot (C\prt D) && ((A\ot B)\ot C)\prt D \\
			&& \text{(P1)} && \\
			& A\ot ((B\ot C)\prt D) && (A\ot (B\ot C))\prt D &
			\arrow["\ao", from=1-1, to=1-3]
			\arrow["\distl", from=1-3, to=1-5]
			\arrow["\id\,\ot\,\distl", from=1-1, to=3-2]
			\arrow["\distl", from=3-2, to=3-4]
			\arrow["\ao\prt\id", from=3-4, to=1-5]
		\end{tikzcd}
	\end{equation*}
	\endgroup
	
	\begingroup
	\renewcommand{\theequation}{\textup{P2}}
	\refstepcounter{equation}\label{eq:P2}
	\begin{equation*}
		\begin{tikzcd}[column sep=-1cm, row sep=0.28cm]
			A\ot ((B\prt C)\ot D) && A\ot (B\prt {(C\ot D)}) && (A\ot B)\prt {(C\ot D)}\\
			&& \text{(P2)} && \\
			& (A\ot (B\prt C))\ot D && ((A\ot B)\prt C)\ot D &
			\arrow["\id\,\ot\,\distr", from=1-1, to=1-3]
			\arrow["\distl", from=1-3, to=1-5]
			\arrow["\ao", from=1-1, to=3-2]
			\arrow["\distl\,\ot\,\id", from=3-2, to=3-4]
			\arrow["\distr", from=3-4, to=1-5]
		\end{tikzcd}
	\end{equation*}
	\endgroup
	
	\begingroup
	\renewcommand{\theequation}{\textup{P3}}
	\refstepcounter{equation}\label{eq:P3}
	\begin{equation*}
		\begin{tikzcd}[column sep=-1cm, row sep=0.28cm]
			A\ot (B\prt {(C\prt D)}) && (A\ot B)\prt {(C\prt D)} && ((A\ot B)\prt C)\prt D \\
			&& \text{(P3)} && \\
			& A\ot ((B\prt C)\prt D) && (A\ot (B\prt C))\prt D &
			\arrow["\distl", from=1-1, to=1-3]
			\arrow["\ap", from=1-3, to=1-5]
			\arrow["\id\,\ot\,\ap", from=1-1, to=3-2]
			\arrow["\distl", from=3-2, to=3-4]
			\arrow["\distl\prt\id", from=3-4, to=1-5]
		\end{tikzcd}
	\end{equation*}
	\endgroup
	
	\begingroup
	\renewcommand{\theequation}{\textup{P4}}
	\refstepcounter{equation}\label{eq:P4}
	\begin{equation*}
		\begin{tikzcd}[column sep=-1cm, row sep=0.28cm]
			(A\prt B)\ot (C\ot D) && A\prt {(B\ot (C\ot D))} && A\prt {((B\ot C)\ot D)} \\
			&& \text{(P4)} && \\
			& ((A\prt B)\ot C)\ot D && (A\prt {(B\ot C)})\ot D &
			\arrow["\distr", from=1-1, to=1-3]
			\arrow["\id\prt\ao", from=1-3, to=1-5]
			\arrow["\ao", from=1-1, to=3-2]
			\arrow["\distr\,\ot\,\id", from=3-2, to=3-4]
			\arrow["\distr", from=3-4, to=1-5]
		\end{tikzcd}
	\end{equation*}
	\endgroup
	
	\begingroup
	\renewcommand{\theequation}{\textup{P5}}
	\refstepcounter{equation}\label{eq:P5}
	\begin{equation*}
		\begin{tikzcd}[column sep=-1cm, row sep=0.28cm]
			(A\prt B)\ot (C\prt D) && ((A\prt B)\ot C)\prt D && (A\prt {(B\ot C)})\prt D \\
			&& \text{(P5)} && \\
			&  A\prt {(B\ot (C\prt D))}  && A\prt {((B\ot C)\prt D)} &
			\arrow["\distl", from=1-1, to=1-3]
			\arrow["\distr\prt\id", from=1-3, to=1-5]
			\arrow["\distr", from=1-1, to=3-2]
			\arrow["\id\prt\distl", from=3-2, to=3-4]
			\arrow["\ap", from=3-4, to=1-5]
		\end{tikzcd}
	\end{equation*}
	\endgroup
	
	\begingroup
	\renewcommand{\theequation}{\textup{P6}}
	\refstepcounter{equation}\label{eq:P6}
	\begin{equation*}
		\begin{tikzcd}[column sep=-1cm, row sep=0.28cm]
			(A\prt {(B\prt C)})\ot D && ((A\prt B)\prt C)\ot D && (A\prt B)\prt {(C\ot D)} \\
			&& \text{(P6)} && \\
			& A\prt {((B\prt C)\ot D)} && A\prt {(B\prt {(C\ot D)})} &
			\arrow["\ap\,\ot\,\id", from=1-1, to=1-3]
			\arrow["\distr", from=1-3, to=1-5]
			\arrow["\distr", from=1-1, to=3-2]
			\arrow["\id\prt\distr", from=3-2, to=3-4]
			\arrow["\ap", from=3-4, to=1-5]
		\end{tikzcd}
	\end{equation*}
	\endgroup

\begin{center}
\begin{tabular}{c c c c}
\Tone &&& \Ttwo \\[2.2cm]
\Tthree &&& \Tfour\\[2.3cm]
\end{tabular}
\end{center}
\end{definition}

\begin{definition}
	A linearly distributive category is \emph{strong} if its distributors are invertible.
\end{definition}

\begin{definition}\label{def:formal-dgm} (Cf. \cite[\S 2.2.1]{DRS26}.)
	A diagram in a linearly distributive category is \emph{formal} if it is obtained from objects and morphisms using composition, the two monoidal products $\otimes$ and $\partimes$, and the associators, unitors, and distributors, without identifying objects or morphisms except via the equalities required by the axioms of linearly distributive categories.
\end{definition}

\section{The \texorpdfstring{$3$}{3}-computad \texorpdfstring{$\highoverline{\ms{LD}}$}{LD}}\label{sec:LD computad}
We fix integers $0\leq k\leq n$ throughout this section. By an \(n\)-category, we mean a strict globular \(n\)-category, which is a standard framework in higher rewriting theory \cite[Def. 14.2.1]{ABGMMM25}. An \(n\)-functor is a morphism of the underlying globular sets that strictly preserves compositions and units. 

Given an $n$-category $X$, its set of $k$-cells is denoted by $X^k$. For $0\leq i<k$, the $i$-source and $i$-target of a $k$-cell $f\in X^k$ are written as $s_i(f),t_i(f)\in X^i$. For two $k$-cells $f,g\in X^k$ with $s_i(f)=t_i(g)$, their $i$-composition is denoted by $f\star_i g\in X^k$. In particular, we compose from right to left. For $k<n$ and $f\in X^k$, we sometimes abusively denote the identity cell $\id_f\in X^{k+1}$ on $f$ by $f\in X^{k+1}$ instead.

\begin{definition}(\cite[\S 1.1.2]{GuMa12}).
	Let \(X\) be an \(n\)-category. A \emph{\(k\)-sphere} of \(X\) is an ordered pair \((x^{k}_1,x^{k}_2)\) of \(k\)-cells in $X$ such that $s_{i}(x^{k}_1)=s_{i}(x^{k}_2)$ and $t_{i}(x^{k}_1)=t_{i}(x^{k}_2)$ for $0\leq i<k$.
\end{definition}
\vspace{-0.41cm}
\begin{figure}[H]
\scalebox{.94}{
\begin{tikzpicture}
\node (A) at (0,0) {$x^{0}_1$};
\node (B) at (3,0) {$x^{0}_2$};
\node (C) at (6,0) {$x^{0}_1$};
\node (D) at (9,0) {$x^{0}_2$};
\draw[->] (C) to [bend left=60] node[midway, above] {$x^{1}_1$} (D);
\draw[->] (C) to [bend right=60] node[midway, below] {$x^{1}_2$} (D);
\node (E) at (12,0) {$x^{0}_1$};
\node (F) at (15,0) {$x^{0}_2$};
\draw[->] (E) to [bend left=60] node[midway, above] {$x^{1}_1$} (F);
\draw[->] (E) to [bend right=60] node[midway, below] {$x^{1}_2$}(F);
 \draw[-{Implies}, double, double distance=2pt] (13.6,.85) to [bend left=60] node[midway, right] {$x^{2}_1$} (13.6,-.85);
\draw[-{Implies}, double, double distance=2pt] (13.4,.85) to [bend right=60] node[midway, left] {$x^{2}_2$}(13.4,-.85);
\end{tikzpicture}
}
\caption{Schematic depictions of a $0$-, a $1$-, and a $2$-sphere.}
\label{fig:spheres}
\end{figure}
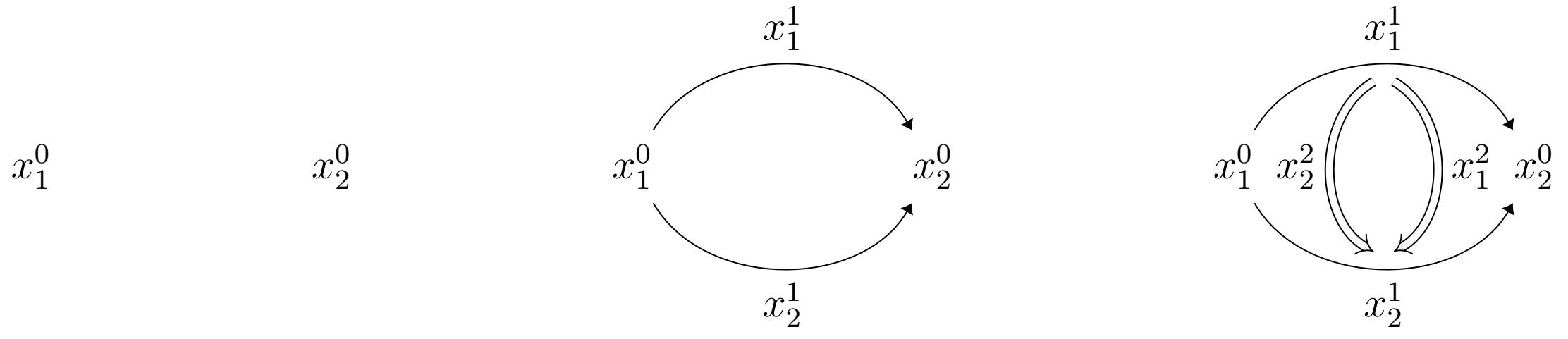

\begin{definition}(\cite[\S 1.1.2]{GuMa12}).
Let \(X\) be an \(n\)-category. A \emph{cellular extension} of \(X\) is a family of \(n\)-spheres of \(X\). The \emph{free \((n+1)\)-category} \(X[\Gamma]\) generated by a cellular extension \(\Gamma\) of \(X\) has, for \(0\leq k\leq n\), the same \(k\)-cells as \(X\), while its \((n+1)\)-cells are freely generated by the \(n\)-spheres \((x_1,x_2)\) of \(\Gamma\), regarded as \((n+1)\)-cells with $n$-source \(x_1\) and $n$-target \(x_2\).
\end{definition}

\begin{remark}
A cellular extension is a \emph{family} and not a \emph{set} of $n$-spheres: one and the same $n$-sphere may occur multiple times in a cellular extension.
\end{remark}

\begin{definition}\label{def:computad} (\cite[\S 1.2]{GuMa12}). The following notions are defined recursively on \(n\in \NNez\):
	\begin{enumerate}[label=(\roman*)]
		\item A \emph{\(0\)-computad} is a set whose elements are called \emph{$0$-cells}. The \emph{free $0$-category} generated by this $0$-computad is the set itself.
		\item An \emph{\((n+1)\)-computad} is a pair \(\Si=(\Sigma_n,\Sigma_{n+1})\), where \(\Sigma_n\) is an \(n\)-computad and \(\Sigma_{n+1}\) is a cellular extension of the free \(n\)-category \(\Si^{\ast}_{n}\) generated by the \(n\)-computad \(\Si_{n}\). The $n$-spheres of \(\Sigma_{n+1}\) are called \emph{$(n+1)$-cells}. The \emph{free \((n+1)\)-category \(\Si^{\ast}\)} generated by \(\Si\) is defined as \(\Si^{\ast}=\Si^{\ast}_n[\Si_{n+1}]\).
	\end{enumerate}
\end{definition}

\begin{remark}
We follow the terminology of \cite{St76,St96,Ba98}; computads are also known as \emph{polygraphs} \cite{B93}.
\end{remark}
For the intuition of the reader, we spell out low-dimensional examples explicitly.
\begin{example}
\begin{enumerate}[label=(\roman*)]
\item A $1$-computad is a rewrite system, in other words, a quiver. Its free $1$-category is the category mentioned in Section \ref{sec:Higher-rewriting} of the introduction.
\item A $2$-computad $(\Si_{1},\Si_{2})$	consists of a quiver $\Si_{1}$ together with a family of $1$-spheres $\Si_{2}$ in the free $1$-category $\Si_{1}^{\ast}$ generated by \(\Si_{1}\).
\item A $3$-computad $(\Si_{2},\Si_{3})$	consists of a $2$-computad $\Si_{2}$ together with a family of $2$-spheres in the free $2$-category $\Si_2^{\ast}$ generated by $\Si_{2}$.
\end{enumerate}
\end{example}

The following \(3\)-computad is central to this paper. We depict its $2$-cells and $3$-cells using string diagrams (read top to bottom, left to right) and surface diagrams, respectively:

\begin{definition}\label{def:3-computad LD}
	The $3$-computad $\ms{LD}$ of \emph{strong linearly distributive categories} has one $0$-cell~\(\ast\), one $1$-cell \(C\) with $\id_C=\twocell{id}$\,, four  $2$-cells $\twocell{ot}$, $\twocell{oteta}$\,,  $\twocell{par}$, $\twocell{pareta}$\,,
	and the following eight $3$-cells:
	\begin{align}\tag{\(\ot\)}
	\scalebox{1}
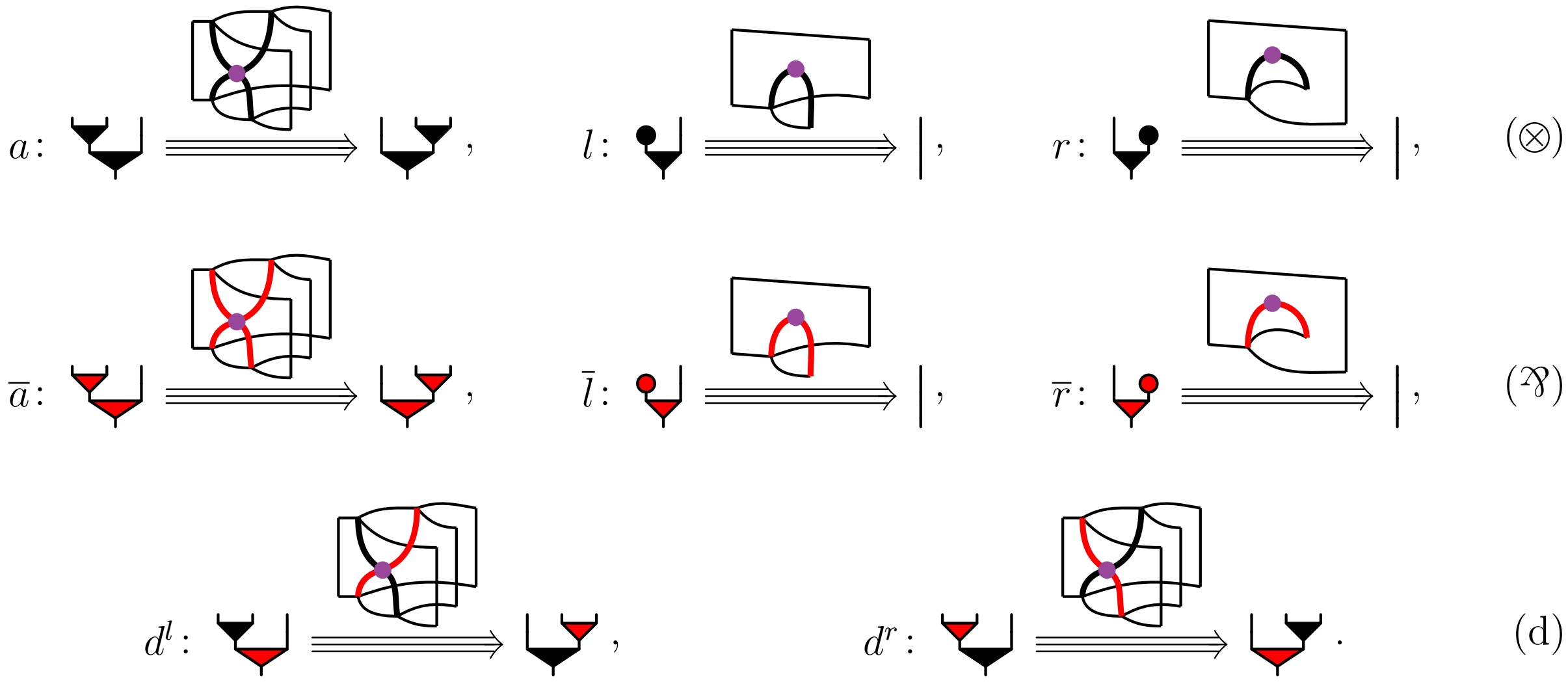
{
	\xymatrix{
		{\aoc\colon\; \twocell{ (ot *0 1) *1 ot }}\;
		\ar@3[rr]^-{\aosurf}
		&& \;{\twocell{ (1 *0 ot) *1 ot } }
	},
	}
	&&
	\scalebox{1}{
	\xymatrix{
	{\loc \colon\; \twocell{ (oteta *0 1) *1 ot }}\;
	\ar@3[rr] ^-{\losurf}
	&& \;{\twocell{id *1 id}}
	},
	}
	&&
	\scalebox{1}{
	\xymatrix{
		{\roc\colon\; \twocell{ (1 *0 oteta) *1 ot }}\;
		\ar@3[rr] ^-{\rosurf}
		&& \;{\twocell{id *1 id}}
	},
	}
	\end{align}
	\begin{align}\tag{\(\prt\)}
		\scalebox{1}{
		\xymatrix{
			{\apc \colon\; \twocell{ (par *0 1) *1 par }}\;
			\ar@3[rr]^-{\apsurf}
			&& \;{\twocell{ (1 *0 par) *1 par } }
		},
		}
		&&
		\scalebox{1}{
		\xymatrix{
			{\lpc\colon\; \twocell{ (pareta *0 1) *1 par }}\;
			\ar@3[rr] ^-{\lpsurf}
			&& \;{\twocell{id *1 id}}
		},
		}
		&&
		\scalebox{1}{
		\xymatrix{\rpc \colon\;
			{\twocell{ (1 *0 pareta) *1 par}}\;
			\ar@3[rr] ^-{\rpsurf}
			&& \;{\twocell{id *1 id}}
		},
		}
	\end{align}
	\begin{align}\label{eq:distri}\tag{d}
		\scalebox{1}{
		\xymatrix{
			{\distlc\colon\;
			 \twocell{ (ot *0 1) *1 par }}\;
			\ar@3[rr] ^-{\distlsurf}
			&& \;{\twocell{ (1 *0 par) *1 ot }}
		},
	}
		&&
		\scalebox{1}{
		\xymatrix{
		{\distrc\colon\;
		\twocell{ (par *0 1) *1 ot }}\;
		\ar@3[rr] ^-{\distrsurf}
		&& \;{\twocell{ (1 *0 ot) *1 par } }
	}.
}	
	\end{align}
\end{definition}

\begin{remark}
The string diagrams $\twocell{ot}$\,,  $\twocell{oteta}$\,, and $s_2(a)$ represent $2$-cells $C \star_0 C \Rightarrow C$, $\id_{\ast} \Rightarrow C$, and $\twocell{ot} \star_1 {\big(\, \twocell{ot} \star_0 \id_C \,\big)}\colon\, C \star_0 C \star_0 C \Rightarrow C$, respectively.
\end{remark}

\begin{remark}
	The $3$-computad $\ms{LD}$ arises from strong linearly distributive categories by orienting their invertible coherence morphisms and omitting their coherence diagrams. When interpreting $\twocell{ot}$ as $\otimes$ and $\twocell{par}$ as $\partimes$, the $3$-cell $d^l$ is oriented oppositely to $\distl$ in \eqref{eq: distributorl}.
\end{remark}

The $3$-computad $\ms{LD}$ does not enjoy all the desired properties; namely, it will turn out to not be confluent in the sense of Definition \ref{def: branch and confl}. Following the Knuth--Bendix completion procedure from \cite[\S 10.2.11]{ABGMMM25}, we thus add the subsequent $3$-cells:

\begin{definition}\label{def: Knuth-Bendix completion LD}
	The $3$-computad $\oLD$, called the \emph{completion} of $\ms{LD}$, has the same $0$-,$1$-,$2$- and $3$-cells as $\ms{LD}$, and additionally the following two $3$-cells:
	\begin{align}\label{eq:units}\tag{compl}
		\scalebox{1}{
		\xymatrix{
			{k^l\colon\;
				\twocell{ (1 *0 oteta *0 1) *1 (1 *0 par) *1 ot }}\;
			\ar@3[rr] ^-{\gammalsurf} 
			&& \;{\twocell{par}}
		},
		}
		&&
		\scalebox{1}{
		\xymatrix{
			{k^r\colon\;
				\twocell{ (1 *0 pareta *0 1) *1 (1 *0 ot) *1 par }}\;
			\ar@3[rr] ^-{\gammarsurf}
			&& \;{ \twocell{ot} }
		}.
		}
	\end{align}
\end{definition}

Proposition \ref{prop:LDCat-algebras} will show that the categorical structures described by a quotient $3$-category built from $\oLD$ are precisely strong linearly distributive categories.

\addtocontents{toc}{\protect\setcounter{tocdepth}{1}}
\section{Computads as higher rewrite systems}

In this section, we recall terminology for computads; readers familiar with higher rewriting may skip this section. Let \(n \geq 0\).

\subsection{Termination}

Termination is defined in terms of a computad's top-dimensional cells:

\begin{definition}\label{def:termination}(\cite[\S 4.1.1]{GuMa}). Let $\Sigma$ be an $(n+1)$-computad.
	\begin{enumerate}[label=(\roman*)]
	\item An $n$-cell \(f\) of the free $(n+1)$-category \(\Si^{\ast}\) generated by $\Sigma$ \emph{rewrites} to an $n$-cell \(g\) of \(\Si^{\ast}\) if there is a non-identity $(n+1)$-cell \(h\) of \(\Sigma^{\ast}\) from \(f\) to \(g\).
	\item A \emph{rewrite sequence} of $\Sigma$ is a finite or infinite sequence \((f_k)\) of \(n\)-cells such that each \(f_k\) rewrites to \(f_{k+1}\).
\end{enumerate}
\end{definition}

\begin{definition} (\cite[\S 4.1.2]{GuMa}).\label{def: termination}
	A computad is \emph{terminating} if any rewrite sequence is finite.
\end{definition}

\subsection{Confluence and convergence}
In this subsection, we single out \emph{convergent} computads. Convergence will allow us to apply Proposition \ref{prop: homotop basis free track} in the proof of our main theorem.
\begin{definition}(Cf. \cite[\S 1.4.1]{GuMa}). 
	Let $\Sigma$ be an $(n+1)$-computad. The \emph{size} \(\|f\|\) of a \(k\)-cell \(f\) in the $(n+1)$-category \(\Sigma^{\ast}\) is the number of \(k\)-cells in \(\Sigma\) that \(f\) is composed of.
\end{definition}
\begin{example}
The $3$-cell $\aoc\star_0 C$ of the free $3$-category \(\ms{LD}^{\ast}\) generated by $\ms{LD}$ has size $1$, where $C$ abusively denotes the identity $3$-cell on the identity $2$-cell on the $1$-cell $C$.
\end{example}

Confluence of computads is defined in terms of their top-dimensional~cells:

\begin{definition}(\cite[\S 4.1.3]{GuMa}).\label{def: branch and confl} Let $\Sigma$ be an $(n+1)$-computad.
	\begin{enumerate}[label=(\roman*)]
	\item A \emph{branching} of \(\Sigma\) is an unordered pair \(\{f,g\}\) of \((n+1)\)-cells of \(\Sigma^{\ast}\) with $s_n(f)=s_n(g)$.
	\item A \emph{branching} \(\{f,g\}\) is \emph{local} if \(\|f\|=\|g\|=1\).
	\item A \emph{confluence} of \(\Sigma\) is an unordered pair $\{f,g\}$ of \((n+1)\)-cells of \(\Sigma^{\ast}\) with $t_n(f)=t_n(g)$.
	\item A branching \(\{f,g\}\) is \emph{confluent} if there exists a confluence \(\{h,k\}\) with \(t_n(f)=s_n(h)\) and \(t_n(g)=s_n(k)\).
	\item The computad \(\Sigma\) is \emph{(locally) confluent} if every (local) branching is confluent.
	\end{enumerate}
\end{definition}

\begin{definition}(\cite[\S 4.1.4]{GuMa}).\label{def: convergent poly}
	A computad is \emph{convergent} if it is terminating and confluent.
\end{definition}

\begin{example}\label{ex: non-confluence}
	In Section~\ref{sec: convergence of compl}, we show that the completion $\highoverline{\ms{LD}}$ is confluent. The $3$-computad $\ms{LD}$ is not confluent, however: two non-confluent local branchings of $\ms{LD}$ are displayed in Figure~\ref{fig: non-confl LD}.
	In fact, these are the only non-confluent \emph{critical} branchings of $\ms{LD}$; see Sections \ref{sec: crit branchs} and \ref{sec: conf}.
	\begin{figure}[H]
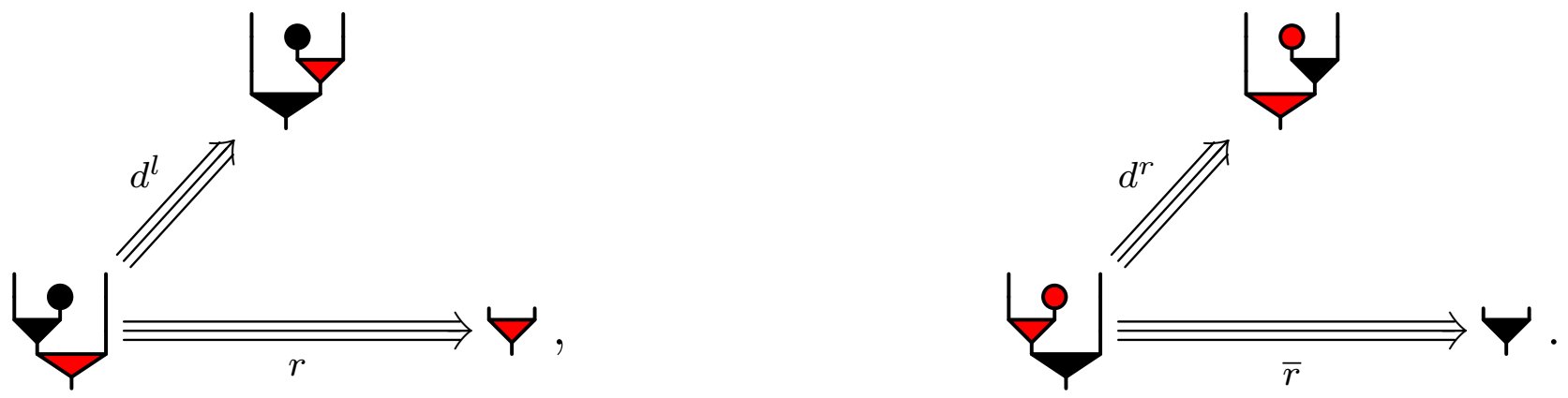

		\centering
		\Critdistlronc \hspace{6em} 
		\Critdistrrpnc
		\caption{The only non-confluent critical branchings of $\ms{LD}$.}
		\label{fig: non-confl LD}
	\end{figure}
\end{example}

\subsection{Critical branchings}\label{sec: crit branchs}
Showing termination is often straightforward. Establishing confluence a priori is harder, since a computad may have infinitely many branchings. The main theorem recalled in this subsection is Theorem~\ref{lem: crit branch lemma}, also known as the \emph{critical branching lemma}, \emph{Newman's lemma}, or the \emph{diamond lemma}. This result reduces the problem of confluence to checking only the critical branchings, of which there are often finitely many. As recalled in Definition~\ref{def: critical branchings}, critical branchings are certain minimal elements with respect to a partial order recalled in Definition~\ref{def: partial order on branchings}. We need the next notion to define this partial order.
\begin{definition}(Cf. \cite[\S 2.1.4]{GuMa}).
Let \(X\) be an \((n+1)\)-category. A \emph{whisker} of an $n$-sphere $x=(x_1,x_2)$ in $X$ is a sequence of pairs $(v_k,w_k)_{1\leq k\leq n}$ of $k$-cells $v_k,w_k$ for which the composite
\begin{equation}\label{eq: whisker}
	v_{n}\star_{n-1}(\,\cdots\, (v_2 \star_{1} (v_1\star_0x\star_0 w_1)\star_{1} w_2) \,\cdots\,) \star_{n-1} w_{n}
\end{equation}
is defined. Here, $x$ is regarded as an $(n+1)$-cell adjoined to $X$, with $n$-source $x_1$ and $n$-target~$x_2$. The composite \eqref{eq: whisker} is viewed as an $(n+1)$-cell in the free $(n+1)$-category \(X^{\leq n}[x]\) generated by the cellular extension $x$ of the truncated $n$-category \(X^{\leq n}\). Note that the expression \eqref{eq: whisker} involves different compositions $\star_i$ of $X$, for $0\leq i\leq n-1$.
\end{definition}

In the following example, we identify each non-negative integer $m$ with the $1$-cell in $\ms{LD}^{\ast}$ obtained by composing $C$ with itself $m$ times along $\star_0$, taking the $0$-fold composite to be \(\id_{\ast}\).

\begin{example}
	A whisker of a $2$-sphere $x=(x_1,x_2)$ in $\ms{LD}^{\ast}$, with $x_1,x_2\colon m \Rightarrow n$, is given by a pair of $2$-cells $v_2\colon p \Rightarrow {k+m+l}$, $w_2\colon {k+n+l} \Rightarrow {q}\,,$ where $k,l,m,n,p,q\geq 0$. This data uniquely determines $1$-cells $v_1,w_1$. For example, the $3$-cell~\eqref{eq: whisker} for a whisker of the $2$-sphere
	\begin{equation}\label{eq:klmn1}
		(x_1=\twocell{ (pareta *0 1) *1 ot }\,,\;\; x_2=\id_1)
	\end{equation}
	is represented diagrammatically in Figure~\ref{fig:whisker}.
\begin{figure}[H]
		\centering
	\whisker
	\caption{The $3$-cell \eqref{eq: whisker} for the $2$-sphere \eqref{eq:klmn1}; $v_k$ and $w_k$ are identity $3$-cells.}
	\label{fig:whisker}
\end{figure}
\end{example}

	Let $W$ be a whisker of an $n$-sphere $x=(x_1,x_2)$ in an \((n+1)\)-category $X$. For any $(n+1)$-cell $f$ in $X$ with $x=(s_n(f),t_n(f))$, write 
	\begin{equation*}
		W[x:=f]
	\end{equation*}
	for the $(n+1)$-cell obtained by substituting $f$ for $x$ in \eqref{eq: whisker}. We now use this substitution to define a partial order on the branchings of a computad.

\begin{definition}(\cite[\S 4.1.5]{GuMa}).\label{def: partial order on branchings}
Let $f=\{f_1,f_2\}$ and $g=\{g_1,g_2\}$ be branchings of an $(n+1)$-computad $\Sigma$. Consider the following $n$-spheres in $\Si^{\ast}$ obtained from the branching $f$:
\begin{equation*}
x\eq (s_{n}(f_1),t_{n}(f_1)) \qquad \text{and} \qquad y\eq (s_{n}(f_2),t_{n}(f_2))\,.
\end{equation*}
We write $f\preccurlyeq g$ if there is a whisker $W$ of both $x$ and $y$ such that 
\begin{equation*}
	W[x:=f_1]\eq g_1 \qquad \text{and} \qquad W[y:=f_2]\eq g_2\,.
\end{equation*}
\end{definition}

\begin{example}
Figure~\ref{fig: partial order} shows two $\preccurlyeq$-comparable branchings of the $3$-computad $\ms{LD}$. The comparison is witnessed by the whisker $v_1\eq w_1 = \id_{0}$, $v_2 = \id_{0} \star_0 \twocell{oteta} \star_0 \id_{0}$, $w_2 \eq \id_{1}$.

\begin{figure}[H]
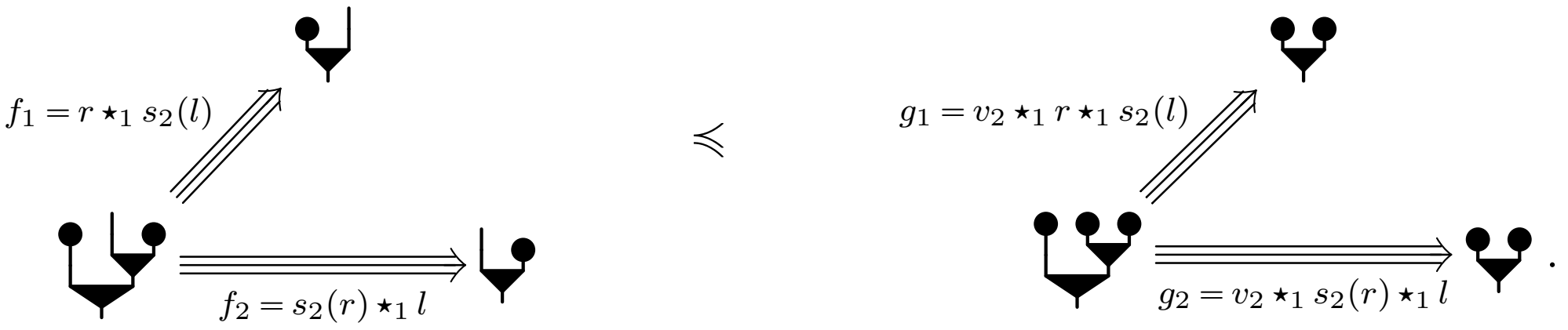

	{\centering
	\begin{tabular}{ccc}
		\Partords 
		& \raisebox{-5ex}{\qquad $\preccurlyeq$\qquad}
		& \Partordb
	\end{tabular}
}
	\caption{Two $\preccurlyeq$-comparable branchings of $\ms{LD}$.}
	\label{fig: partial order}
\end{figure}
\end{example}

\begin{definition}(Cf. \cite[\S 4.1.5]{GuMa}).\label{def: critical branchings}
	Let $\Si$ be an $(n+1)$-computad. A \emph{minimal branching} of $\Si$ is a minimal element with respect to $\preccurlyeq$. A \emph{critical branching} of $\Si$ is a minimal local branching that is neither of the form $\{f,f\}$, for an $(n+1)$-cell $f$ of $\Sigma$, nor of the form $\{f\star_i s_n(g),s_n(f)\star_i g\}$, for $(n+1)$-cells $f,g$ of $\Sigma$ and $0\leq i \leq n-1$.
\end{definition}

The main result recalled in this subsection is the \emph{critical branching lemma}:

\begin{theorem}\label{lem: crit branch lemma} \emph{(\cite[Thm. 2.5.4]{Gu19}).}
	A terminating $(n+1)$-computad \(\Sigma\) is confluent if and only if every critical branching of \(\Sigma\) is confluent.
\end{theorem}

\subsection{Aspherical higher track categories}\label{sec: hombasis}

For a cellular extension  $\Gamma$ of an $n$-category $X$, recall the following terminology from \cite[\S 1.1.2]{GuMa12}.

\begin{definition}
The \emph{quotient} $X/\Gamma$ is the $n$-category obtained from $X$ by identifying, for all $n$-spheres 
$(x_1,x_2)\in \Gamma$, the $n$-cells $x_1$ and 
$x_2$.
\end{definition}

\begin{definition}
The \emph{free track $(n+1)$-category} generated by $\Gamma$ over $X$ is the quotient
\begin{equation*}
	X(\Gamma) \eq X[\Gamma,\Gamma^{-}]\,/\,{\operatorname{Inv}(\Gamma)}\,,
\end{equation*}
where $\Gamma^{-}$ denotes the following cellular extension of $X$
\begin{align*}
	\Gamma^{-} &\eq \big\{\gamma^{-}=(x_2,x_1) \;\;\vert\;\; \gamma=(x_1,x_2)\in \Gamma\big\}\,,
\end{align*}
$X[\Gamma,\Gamma^{-}]$ is the free $(n+1)$-category generated by the cellular extension $\Gamma \cup \Gamma^{-}$ of $X$, and~$\operatorname{Inv}(\Gamma)$~denotes the following cellular extension of $X[\Gamma,\Gamma^{-}]$
\begin{align*}
	\operatorname{Inv}(\Gamma) &\eq \big\{(\gamma \star_n \gamma^{-},\id_{s_{n-1}(\gamma)})\,,\; (\gamma^{-} \star_n \gamma,\id_{t_{n-1}(\gamma)}) \;\;\vert\;\; \gamma \in \Gamma\big\}\,.
\end{align*}
\end{definition}

\begin{remark}
	By construction, the top cells of the $(n+1)$-category $X(\Gamma)$ are invertible.
\end{remark}

Let $\Sigma$ be a $(n+1)$-computad, and let $\Sigma^\top=\Sigma_n^{\ast}(\Sigma_{n+1})$ be the free track {$(n+1)$-category} generated by $\Sigma_{n+1}$. Assume that $\Sigma$ is convergent. For each critical branching $\{f,g\}$ in $\Sigma$, choose a confluence $\{h,k\}$. Let $\Gamma_\Sigma$ be the set of the resulting $(n+1)$-spheres $(h\star_n f,k\star_n g)$.
\begin{proposition}\emph{(\cite[Prop. 4.3.4]{GuMa}).}\label{prop: homotop basis free track}
	For any such choice of confluences, the quotient $(n+1)$-category $\Sigma^\top/\,\Gamma_\Sigma$ is \emph{aspherical}: any two parallel $(n+1)$-cells are equal. One calls $\Gamma_\Sigma$ a \emph{homotopy basis} of the track $(n+1)$-category $\Sigma^\top$.
\end{proposition}
 The proof of \cite[Prop. 4.3.4]{GuMa} relies on the invertibility of the top-dimensional cells in~$\Sigma^\top$. This is one reason why we restrict to \emph{strong} linearly distributive categories.

\addtocontents{toc}{\protect\setcounter{tocdepth}{2}}
\section{Convergence of \texorpdfstring{$\highoverline{\ms{LD}}$}{LD}}\label{sec: convergence of compl}

In this section, we show that the $3$-computad 
$\oLD$ of Definition~\ref{def:3-computad LD} is convergent. By Theorem~\ref{lem: crit branch lemma}, it suffices to prove that $\oLD$ is terminating and that all its critical branchings are confluent. We now establish these properties.

\subsection{Termination}\label{subs:termination}
Let \(\oLD_{2}\) denote the truncation of \(\oLD\) to a \(2\)-computad, and let \(\oLD_{2}^{\ast}\) be the free \(2\)-category it generates. We identify a non-negative integer $m$ with the $1$-cell given by the $m$-fold \(0\)-composition of the generating $1$-cell \(C\), where the $0$-fold composition is \(\id_{\ast}\).

Let \(\Ord\) be the one-object \(2\)-category of posets and monotone functions. Its \(0\)-composition is the cartesian product \(\times\) of posets. Its \(1\)-composition is function composition \(\circ\).

\begin{definition} (Cf. \cite[\S 2.4]{GuMa}).
Let \(X\colon \oLD_{2}^{\ast} \to \Ord\) be a \(2\)-functor. An \emph{$X$-valued derivation} of \(\oLD_{2}^{\ast}\) is an assignment that sends each \(2\)-cell \(F\colon m \Rightarrow n\) of \(\oLD_{2}^{\ast}\) to a monotone~function
\begin{equation*}
	d(F)\colon \; X(C)^{\times m} \;\longrightarrow\; \ZZ\,,
\end{equation*}
subject to the following compatibility conditions with respect to the compositions $\star_0$ and $\star_1$:

\begin{enumerate}
	\item for \(2\)-cells \(F\colon m \Rightarrow n\) and \(G\colon k \Rightarrow l\), we have
	\begin{equation}\label{eq: derivation composition 0}
		d(F \star_0 G){(i_1,\ldots,i_{m+k})}
		\eq
		d(F){(i_1,\ldots, i_m)} + d(G){(i_{m+1},\ldots, i_{m+k})}\,.
	\end{equation}
	\item for \(2\)-cells \(G\colon n \Rightarrow k\) and \(F\colon k \Rightarrow m\), we have
	\begin{equation}\label{eq: derivation composition 1}
		d(F \star_1 G){(i_1,\ldots,i_{n})} 
		\eq
		(d(F)\circ X(G)){(i_1,\ldots,i_n)} + d(G){(i_1,\ldots, i_n)}\,.
	\end{equation}
\end{enumerate}
\end{definition}

\begin{remark}\label{rem: d(id)=0}
Condition \eqref{eq: derivation composition 1} implies that $d(\id_n)=0$ for every $1$-cell $n$ of \(\oLD_{2}^{\ast}\). Therefore, $d$ respects the interchange law. That is, for all $2$-cells $F\colon k \Rightarrow l$ and $G\colon m \Rightarrow n$ of \(\oLD_{2}^{\ast}\), we have
\begin{equation*}
	d\big((\id_l \star_0 G) \star_1 (F\star_0 \id_m)\big)\eq 
	d\big(F \star_0 G\big) \eq
	d\big((F \star_0 \id_n) \star_1 (\id_k\star_0 G)\big)
	\,.
\end{equation*}
\end{remark}

More generally, derivations are defined for any $n$-category; see \cite[\S 2.4]{GuMa}. In this setting:

\begin{proposition}\label{prop:termination lemma} \emph{(Cf. \cite[Thm. 4.2.1]{GuMa}).}
	Let $\Sigma$ be a $3$-computad such that there exist:
	\begin{itemize}
		\item a $2$-functor $X:\Sigma^{\ast}_2 \to \Ord$ such that, for any $1$-cell $f$ in $\Si$, the underlying set of $X(f)$ is non-empty and, for any $3$-cell $g$ in $\Si$, the inequality $X(s_2(g))\geq X(t_2(g))$ holds\,,
		\item and an $X$-valued derivation $d$ of $\Sigma^{\ast}_2$ such that, for any $2$-cell $F$ in $\Sigma^{\ast}_2$, we have $d(F)\geq 0$ and, for any $3$-cell $g$ in $\Si$, the strict inequality $d(s_2(g))>d(t_2(g))$ holds\,.
	\end{itemize}
	
	\noindent Then the $3$-computad $\Si$ is terminating.
\end{proposition}

\pagebreak

\begin{lemma}\label{lemma:termination oLD}
The \(3\)-computad \(\oLD\) is terminating.
\end{lemma}
\begin{proof}
Define the \(2\)-functor \(X\colon \oLD_{2}^{\ast} \to \Ord\) on the generating cells of $\oLD_{2}^{\ast}$ by
	\begin{align*}
	&X\big(\,C\,\big) \eq \NN\setminus\{0\}\,,
	\qquad &X\big(\,\twocell{ot}\,\big){(i,j)} \eq X\big(\,\twocell{par}\,\big){(i,j)} \eq i+j \,,
	\qquad &X\big(\,\twocell{oteta}\,\big) \eq X\big(\,\twocell{pareta}\,\big) \eq 1 \,.
\end{align*}
The following assignment extends uniquely to an $X$-valued derivation \(d\) of \(\wSi\):
\begin{align*}
	d\big(\,\twocell{ot}\,\big){(i,j)} \eq d\big(\,\twocell{par}\,\big){(i,j)}\eq i \,, 
	\qquad d\big(\,\twocell{oteta}\,\big) \eq d\big(\,\twocell{pareta}\,\big) \eq 0 \,.
\end{align*}
By definition and by Equations \eqref{eq: derivation composition 0}–\eqref{eq: derivation composition 1}, $d(F)\geq 0$ for every $2$-cell $F$ in $\oLD^{\ast}_2$. We verify~that
\begin{equation*}\label{eq:ddecreasing}
	X(s_2(g))\geq X(t_2(g)) \quad \text{and} \quad d(s_2(g))> d(t_2(g))
\end{equation*}
for any of the ten \(3\)-cells \(g\) of \(\oLD\). First, we consider \(X\). For $\aoc,\loc,\roc,\apc,\lpc,\rpc$, we have:
\begin{gather*}
	X\Big(\,\twocell{ (ot *0 1) *1 ot }\,\Big){(i,j,k)} 
	\eq 
	X\Big(\,\twocell{ (par *0 1) *1 par }\,\Big){(i,j,k)} 
	\eq
	i+j+k
	\eq
	 X\Big(\,\twocell{ (1 *0 par) *1 par }\,\Big){(i,j,k)} 
	\eq
	X\Big(\,\twocell{ (1 *0 ot) *1 ot }\,\Big){(i,j,k)}\,,\\
	X\Big(\,{\twocell{ (oteta *0 1) *1 ot }}\,\Big){(i)} 
	\eq
	X\Big(\,{\twocell{ (1 *0 oteta) *1 ot }}\,\Big){(i)}
	\eq X\Big(\,{\twocell{ (pareta *0 1) *1 par }}\,\Big){(i)} 
	\eq X\Big(\,{\twocell{ (1 *0 pareta) *1 par}}\,\Big){(i)}  
	\eq i+1 
	\,>\, i
	\eq X\Big(\:{\twocell{id *1 id}}\:\Big){(i)}\,.
\end{gather*}
For the other $3$-cells $\distlc, \distrc$ and $k^l,k^r$, we compute:
\begin{gather*}
	X\Big(\,\twocell{ (ot *0 1) *1 par }\,\Big){(i,j,k)} 
	\eq 
	X\Big(\,\twocell{ (par *0 1) *1 ot }\,\Big){(i,j,k)} 
	\eq
	i+j+k
	\eq
	X\Big(\,\twocell{ (1 *0 par) *1 ot }\,\Big){(i,j,k)} 
	\eq
	X\Big(\,\twocell{ (1 *0 ot) *1 par }\,\Big){(i,j,k)}\,,\\
	X\Big(\,\twocell{ (1 *0 oteta *0 1) *1 (1 *0 par) *1 ot }\,\Big){(i,j)} 
	\eq 
	X\Big(\,\twocell{ (1 *0 pareta *0 1) *1 (1 *0 ot) *1 par }\,\Big){(i,j)} 
	\eq
	i+1+j
	\,>\,
	i+j
	\eq 
	X\Big(\,\twocell{par}\,\Big){(i,j)}
	\eq
	X\Big(\,\twocell{ot}\,\Big){(i,j)}\,.
\end{gather*}
Second, we consider $d$. Using Remark \ref{rem: d(id)=0}, we find for $\aoc,\loc,\roc,\apc,\lpc,\rpc$, and for all \(i,j,k\in \NN\setminus\{0\}\):
{\centering
	\makebox[\textwidth][c]{
	\resizebox{0.99\textwidth}{!}{
	\begin{minipage}{\textwidth}
	\begin{gather*}
	d\Big(\,\twocell{ (ot *0 1) *1 ot }\,\Big){(i,j,k)} 
	\eq 
	d\Big(\,\twocell{ (par *0 1) *1 par }\,\Big){(i,j,k)} 
	\eq
	2i+j
	\,>\, i+j
	\eq
	d\Big(\,\twocell{ (1 *0 par) *1 par }\,\Big){(i,j,k)}
	\eq d\Big(\,\twocell{ (1 *0 ot) *1 ot }\,\Big){(i,j,k)}\,,\\
	d\Big(\,{\twocell{ (1 *0 oteta) *1 ot }}\,\Big){(i)} 
	\eq
	d\Big(\,{\twocell{ (1 *0 pareta) *1 par}}\,\Big){(i)} 
	\eq i
	\,>\, 0
	\eq d\Big(\:{\twocell{id *1 id}}\:\Big){(i)}\,,\\
	d\Big(\,{\twocell{ (oteta *0 1) *1 ot }}\,\Big){(i)} 
	\eq
	d\Big(\,{\twocell{ (pareta *0 1) *1 par }}\,\Big){(i)} 
	\eq 1
	\,>\, 0
	\eq d\Big(\:{\twocell{id *1 id}}\:\Big){(i)}\,.
	\end{gather*}
\end{minipage}
}
}
}

\noindent Finally, for the other $3$-cells $\distlc, \distrc$ and $k^l,k^r$, we have:

{\centering
	\makebox[\textwidth][c]{
		\resizebox{0.99\textwidth}{!}{
			\begin{minipage}{\textwidth}
\begin{gather*}
	d\Big(\,\twocell{ (ot *0 1) *1 par }\,\Big){(i,j,k)} 
	\eq 
	d\Big(\,\twocell{ (par *0 1) *1 ot }\,\Big){(i,j,k)} 
	\eq
	2i+j
	\,>\, i+j
	\eq
	d\Big(\,\twocell{ (1 *0 par) *1 ot }\,\Big){(i,j,k)}
	\eq d\Big(\,\twocell{ (1 *0 ot) *1 par }\,\Big){(i,j,k)}\,,\\
	d\Big(\,\twocell{ (1 *0 oteta *0 1) *1 (1 *0 par) *1 ot }\,\Big){(i,j)} 
	\eq 
	d\Big(\,\twocell{ (1 *0 pareta *0 1) *1 (1 *0 ot) *1 par }\,\Big){(i,j)} 
	\eq
	i+1
	\,>\,
	i
	\eq 
	d\Big(\,\twocell{par}\,\Big){(i,j)}
	\eq
	d\Big(\,\twocell{ot}\,\Big){(i,j)}\,.
\end{gather*}
\end{minipage}
}
}
}
 
The termination of \(\oLD\) now follows from Proposition \ref{prop:termination lemma}.
\end{proof}

\subsection{Confluence}\label{sec: confluence}\label{sec: conf}

A combinatorial analysis of the local branchings of $\oLD$ identifies exactly $36$~critical branchings. We verify that each is confluent by giving an explicit confluence, shaded grey in the figures. Identity $2$-cells are omitted from composites; for instance, in the first critical branching, we write $\distlc$ instead of $(\twocell{oteta}\star_0 C\star_0 C)\star_1\distlc$.

\smallskip

We begin with the confluent critical branchings of the uncompleted computad $\ms{LD}$. Those involving both black and red $2$-cells, representing $\otimes$-structures and $\partimes$-structures, are:
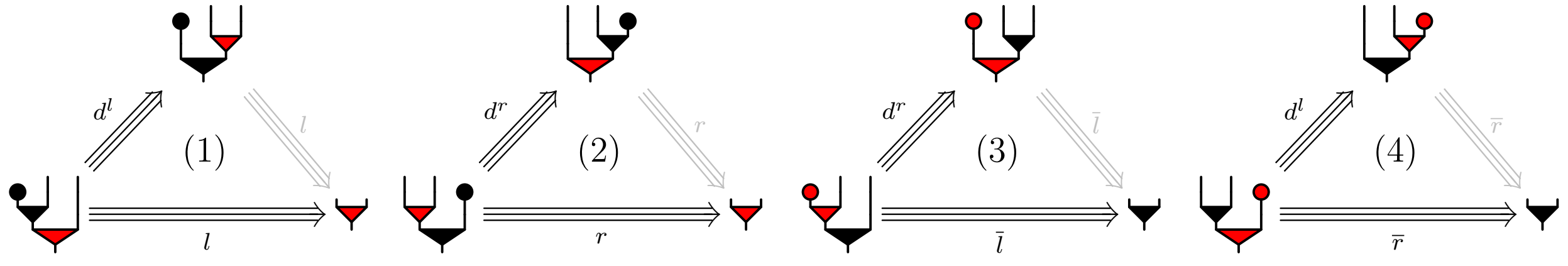
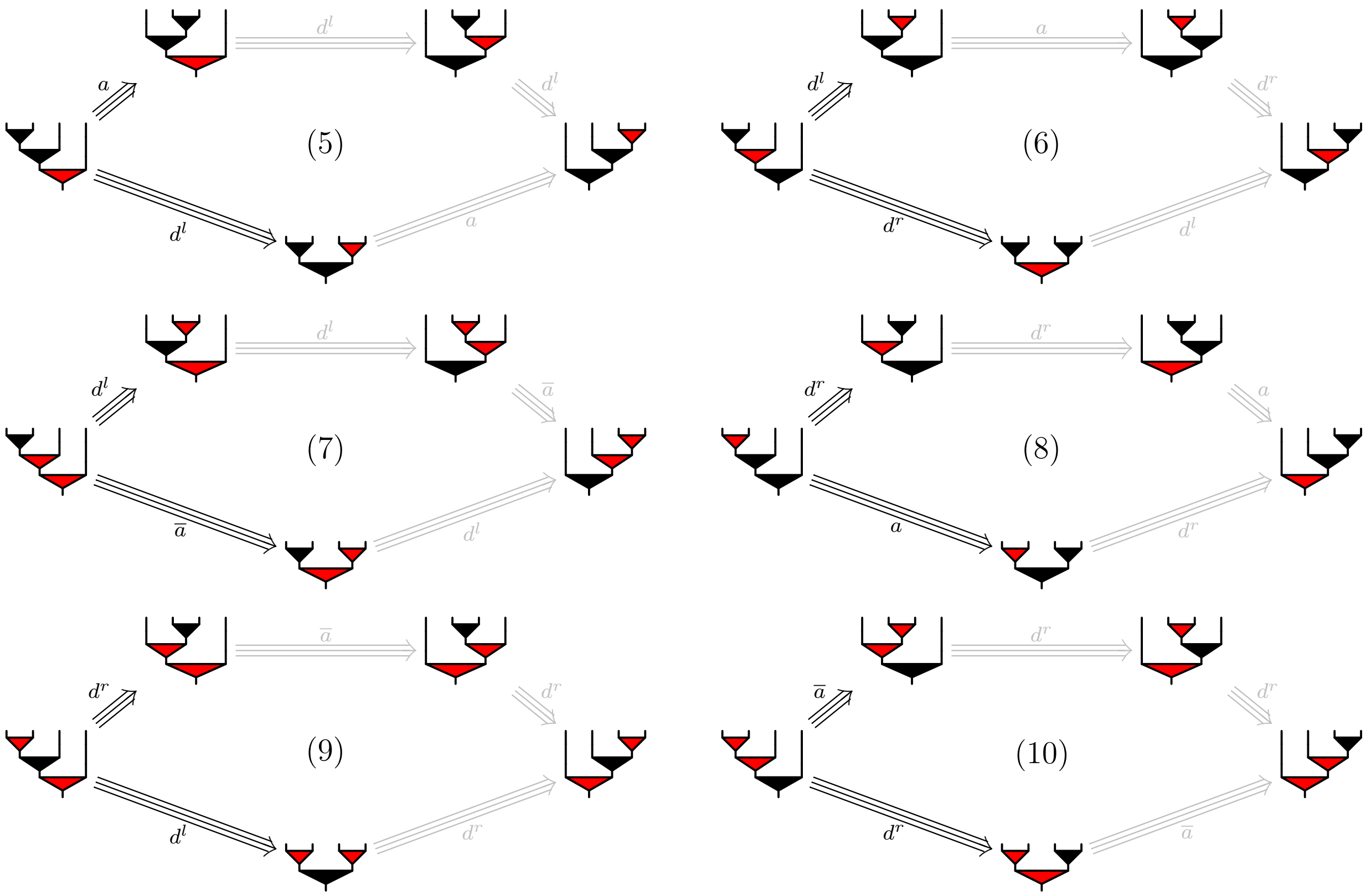
\vspace{0.1cm}

{\centering
\hspace{-0.3cm}
\begin{tabularx}{\textwidth}{ m{3.8cm} m{3.84cm} m{3.84cm} m{3.84cm} }
\Critdistllo & \Critdistrro & \Critdistrlp & \Critdistlrp\\[2.5 cm]
\end{tabularx}
}

{\centering
\begin{tabularx}{\textwidth}{ m{8.2cm} m{8.2cm} }
	\Critdistlao & \Critaodistl \\[3 cm]
\end{tabularx}
}

{\centering
	\begin{tabularx}{\textwidth}{ m{8.2cm} m{8.2cm} }
		\Critdistlap & \Critdistrao \\[3 cm]
	\end{tabularx}
}
{\centering
	\begin{tabularx}{\textwidth}{ m{8.2cm} m{8.2cm} }
		\Critdistldistr & \Critapdistr\\[3 cm]
	\end{tabularx}
}

\vspace{0.2cm}

The critical branchings of $\ms{LD}$ involving only the black \(2\)-cells are:

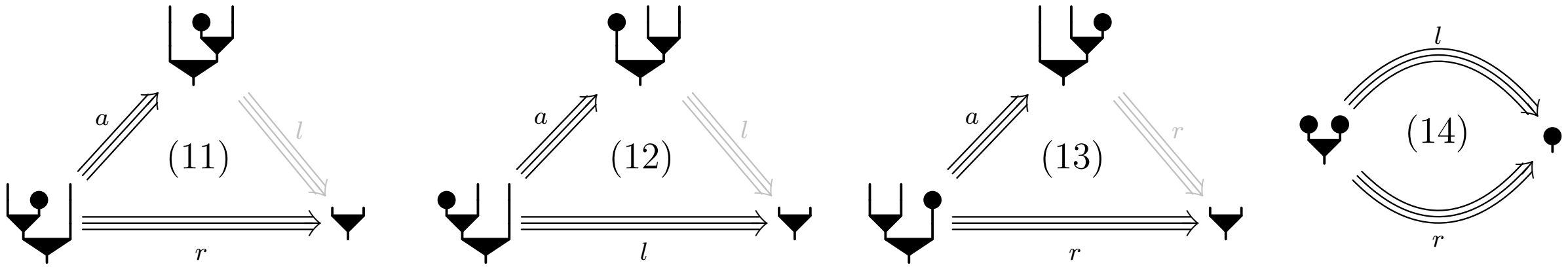
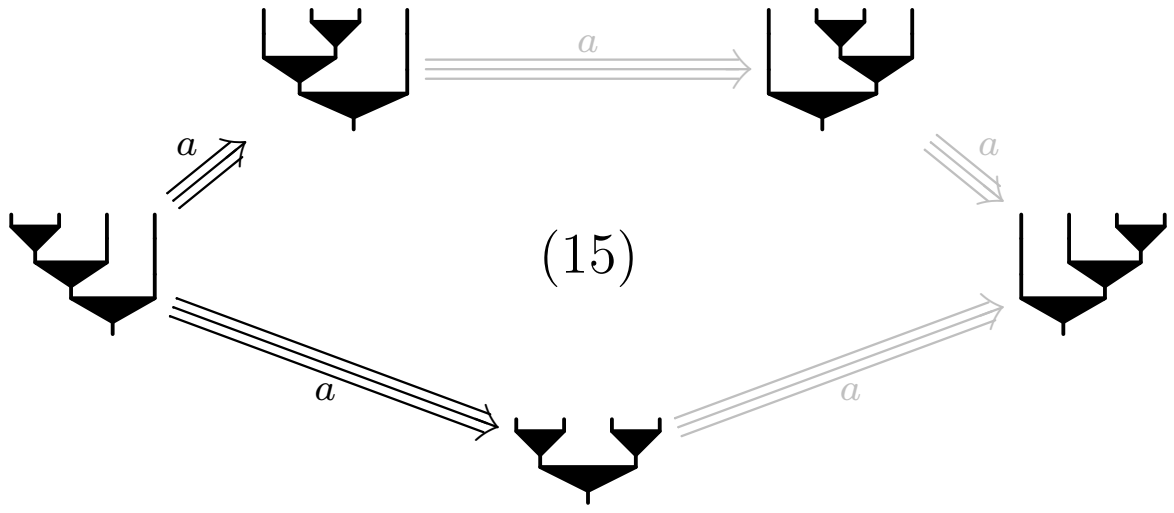
\vspace{0.3cm}

{\centering
\begin{tabularx}{\textwidth}{ m{4cm} m{4cm} m{4cm} m{4cm} }
	\Critaoro &\Critaolo & \Critaororo & \raisebox{-10mm}{\Critloro}
\end{tabularx}
\vspace{0cm}
\begin{longtable}{ c }
\Critaoao
\end{longtable}
}

The critical branchings of $\ms{LD}$ involving only the red \(2\)-cells are:

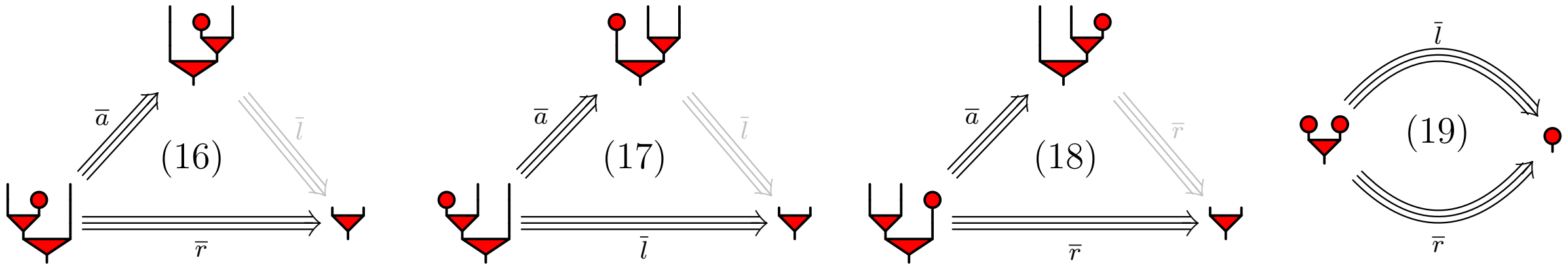
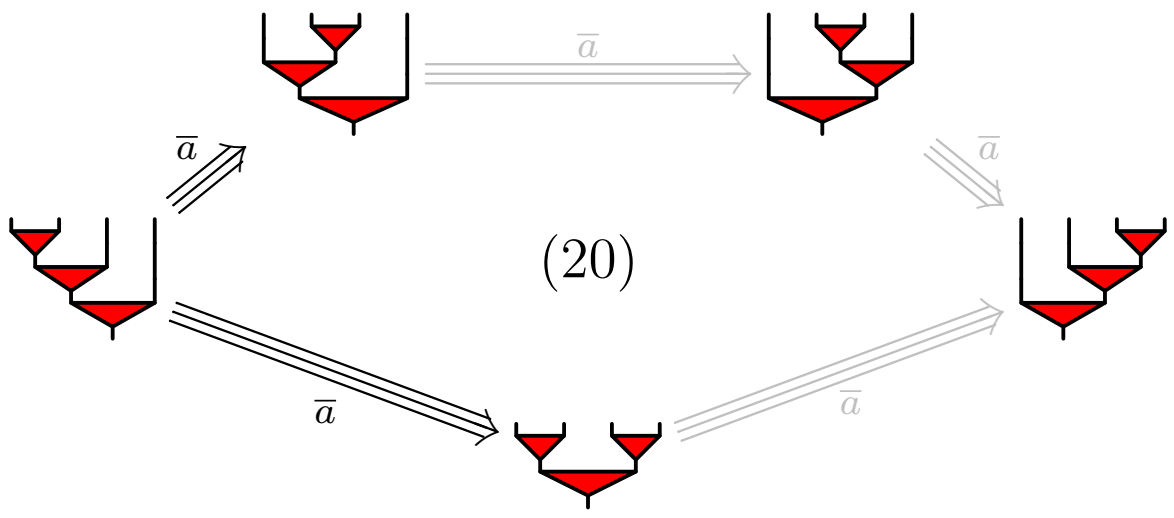
\vspace{0.4cm}

{\centering
\begin{tabularx}{\textwidth}{ m{4cm} m{4cm} m{4cm} m{4cm} }
		\Critaprp &\Critaplp & \Critaprprp & \raisebox{-10mm}{\Critlprp}
\end{tabularx}
\vspace{-0.2cm}
\begin{longtable}{ c }
		\Critapap
	\end{longtable}
}

\vspace{-0.4cm}

Second, there are the non-confluent critical branchings of $\ms{LD}$ from Example \ref{ex: non-confluence} which become confluent in the completion $\oLD$:

\vspace{-0.6cm}

{\centering
\begin{longtable}{ @{}l@{\hspace{3.5cm}}r@{} }
	\Critdistlroncr & \Critdistrrpncr
\end{longtable}
}

\vspace{-0.5cm}

In the terminology of \cite[\S 5.1.1]{GuMa}, these two critical branchings are neither inclusion nor regular, left-indexed, right-indexed, or multi-indexed.

\smallskip

Third, there are the following critical branchings involving $k^l$ but not $k^r$:

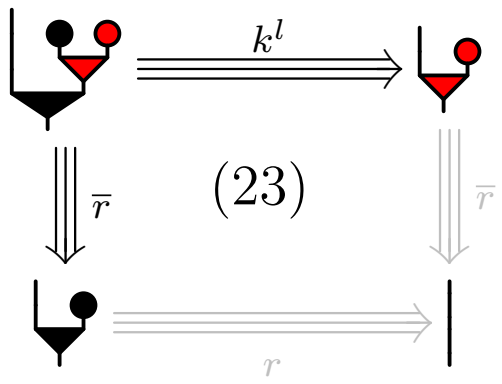
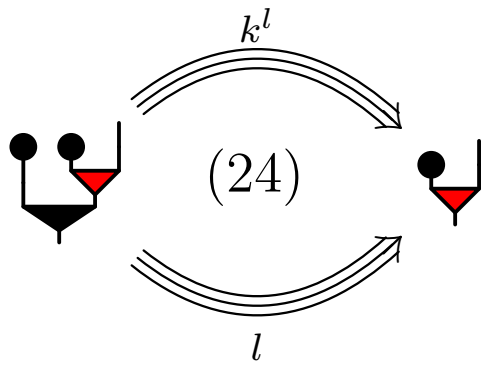
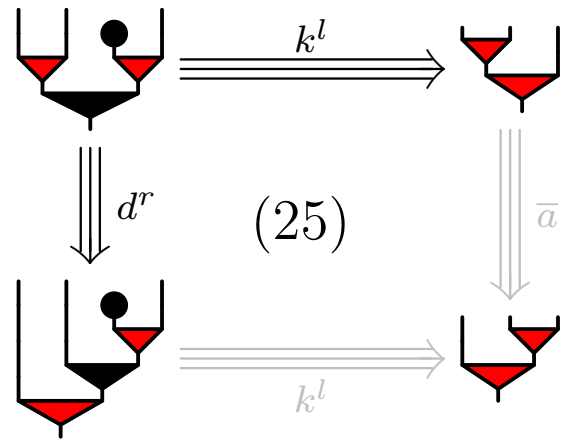
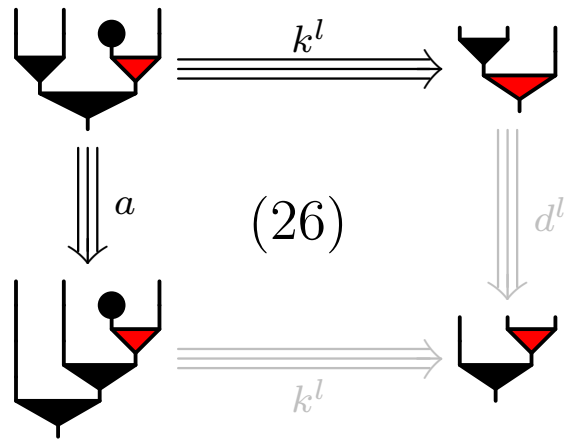
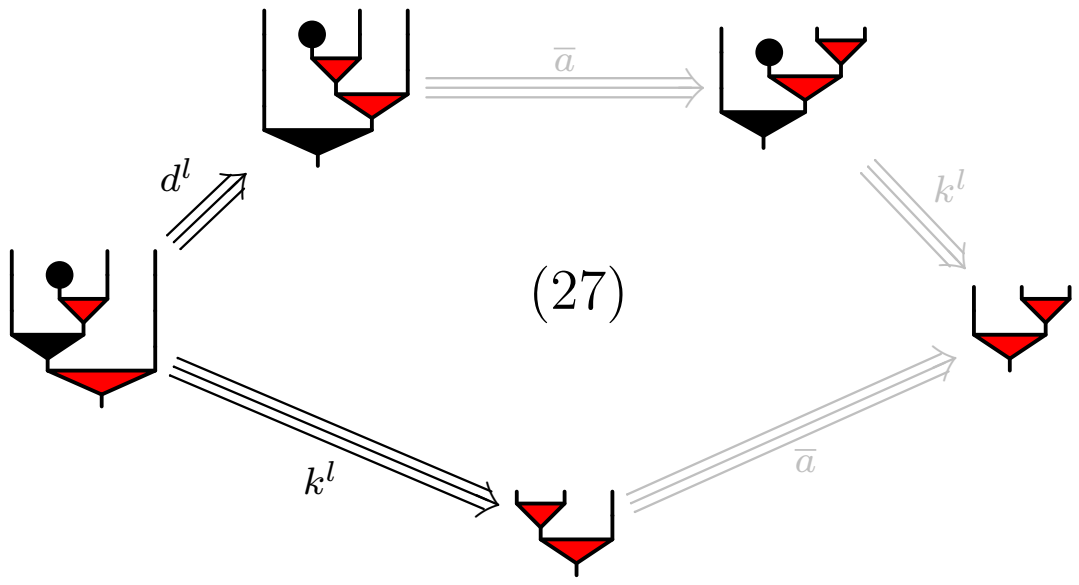
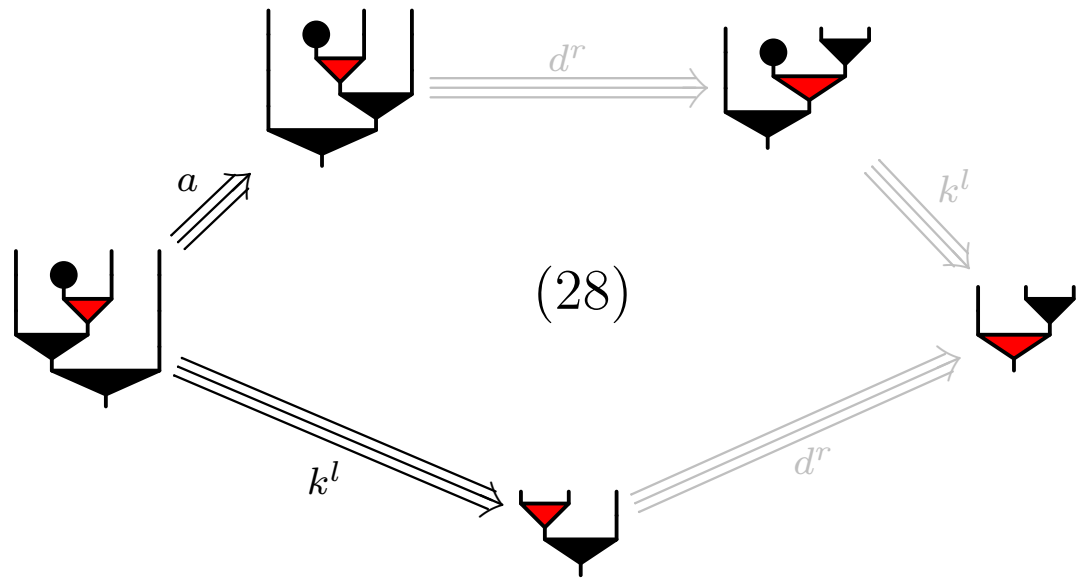
\vspace{0.4cm}

{\centering
\begin{tabularx}{\textwidth}{ m{3.9cm} m{3.9cm} m{3.9cm} m{3.9cm} }
	\Critglrp & \hspace{-10pt}\raisebox{-10mm}{\Critgllo} & \hspace{-15pt}\Critgldistr &\Critglao\\[2.7 cm]
\end{tabularx}
}
{\centering
	\begin{tabularx}{\textwidth}{ m{8.5cm} m{8.5cm} }
			\Critgldistl & \CritglaoP
	\end{tabularx}
}

\vspace{0.5cm}

Interchanging the colours gives the following critical branchings involving $k^r$ but not $k^l$:

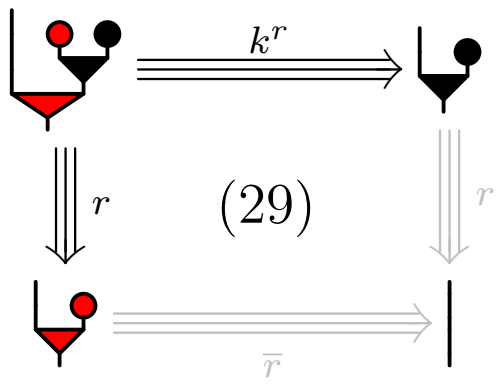
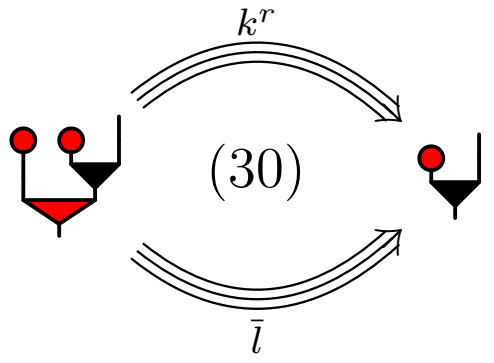
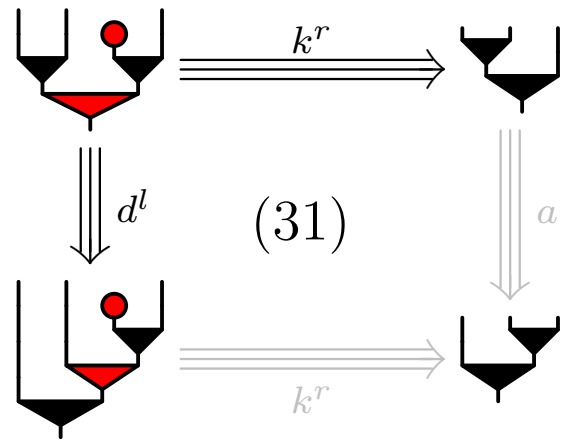
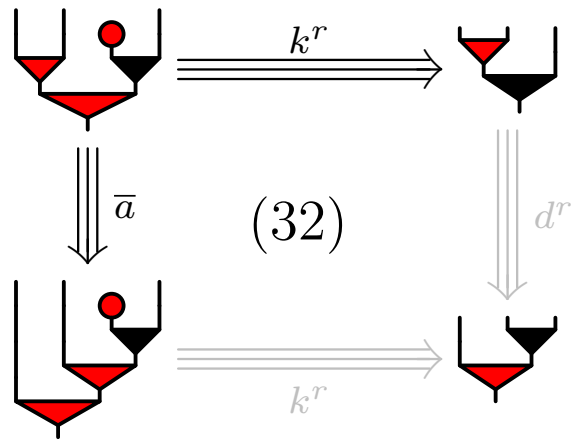
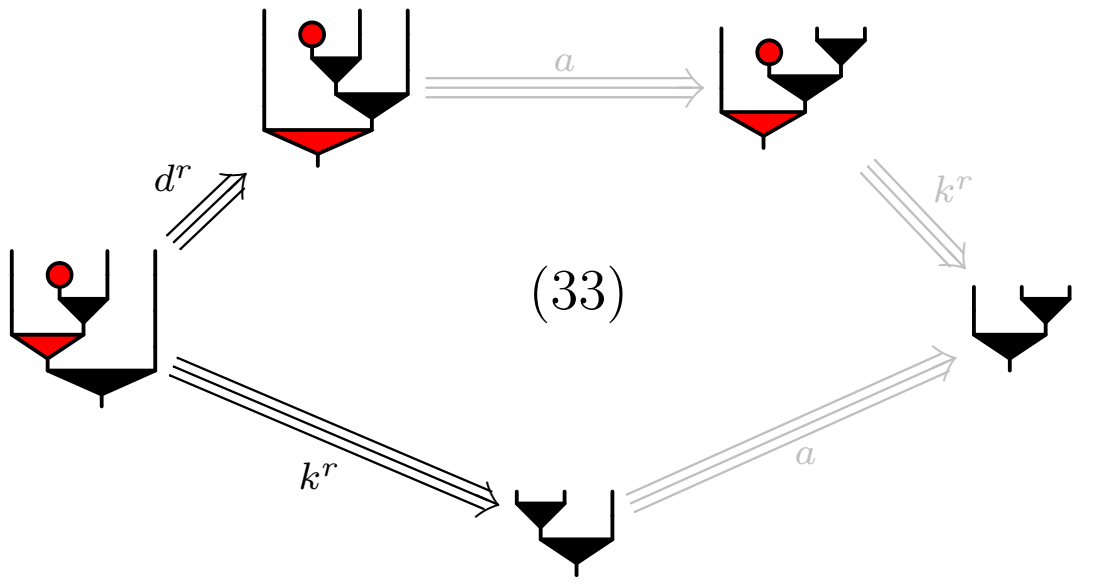
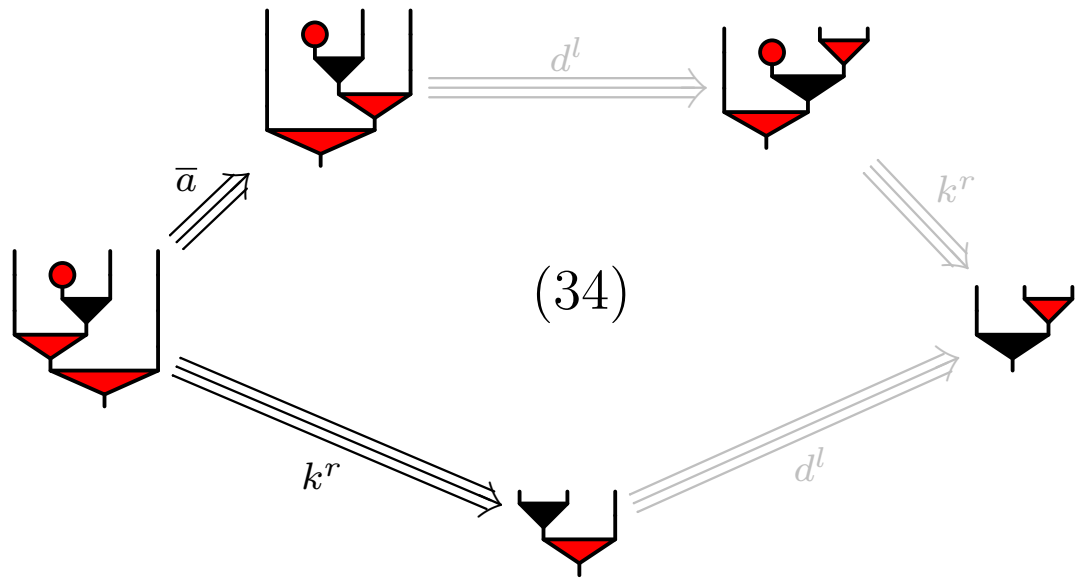
\vspace{0.5cm}

{\centering
\begin{tabularx}{\textwidth}{ m{3.9cm} m{3.9cm} m{3.9cm} m{3.9cm} }
	\Critgrro & \hspace{-10pt}\raisebox{-10mm}{\Critgrlp} & \hspace{-15pt}\Critgrdistl & \Critgrap
\end{tabularx}
}

{\centering
	\begin{tabularx}{\textwidth}{ m{8.2cm} m{8.2cm} }
	\Critgrdistr & \CritgrapP \\[3 cm]
\end{tabularx}
}

\vspace{0.4cm}

Fourth, there are the following critical branchings involving both $k^l$ and $k^r$:

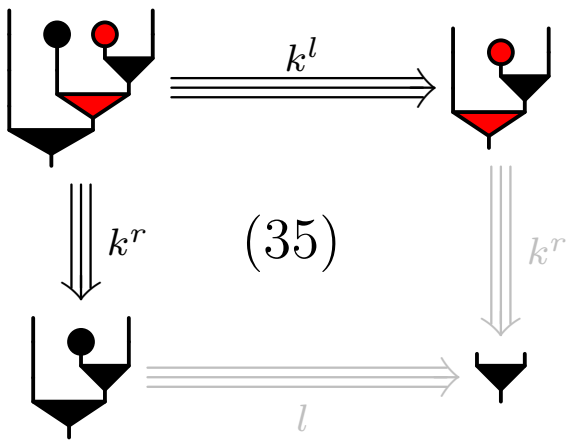
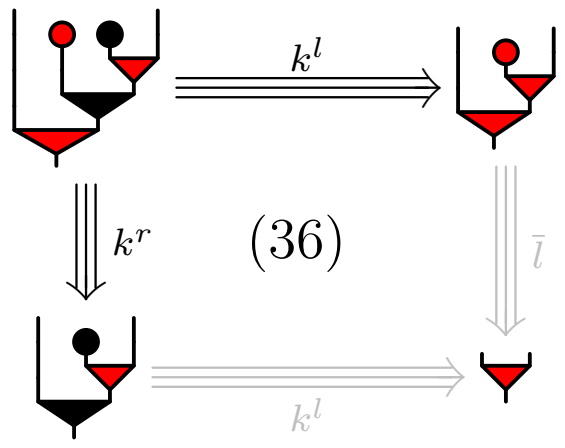
\begin{longtable}{ @{}l@{\hspace{3.5cm}}r@{} } 
	\raggedleft\Critglgro & \Critglgrt
\end{longtable}

The confluence of the critical branchings, together with Theorem \ref{lem: crit branch lemma} and Lemma \ref{lemma:termination oLD} imply:
\begin{proposition}\label{prop: completed computad converg}
The completed \(3\)-computad \(\oLD\) is convergent.
\end{proposition}

\section{Strong linearly distributive coherence}

Following Section~\ref{sec: hombasis}, let $\oLD^{\top}$ be the free track $3$-category generated by the $3$-cells of $\oLD$. Let $\Gamma_{\oLD}$ denote the set consisting of the $3$-spheres in $\oLD$ arising from all its $36$ critical branchings from Section~\ref{sec: confluence}. Let $\mathbf{LDCat}$ denote the quotient $3$-category $\oLD^{\top}\!/\,\Gamma_{\oLD}$.

We use $3$-categories to encode categories equipped with structure:

\begin{definition}\label{def:algebra over tcategory}(Cf. \cite[\S 0]{GuMa12}).
	Let \(X\) be a \(3\)-category. An \emph{$X$-algebra} is a $3$-functor $X \to \mathbf{Cat}$, where $\mathbf{Cat}$ is the one-object \(3\)-category whose $1$-, $2$-, and $3$-cells are categories, functors, and natural transformations, respectively. Its compositions are given by cartesian product, horizontal composition, and vertical composition. 
\end{definition}

The $3$-computad $\oLD$ encodes strong linearly distributive categories:

\begin{proposition}\label{prop:LDCat-algebras}
	Strong linearly distributive categories are in bijection with $\mathbf{LDCat}$-algebras. 
\end{proposition}
\begin{proof}
	The $\mathbf{LDCat}$-algebra $F$ corresponding to a strong linearly distributive category $\cC$, with monoidal structures \((\otimes,1)\) and \((\prt,K)\), is defined on the generators of $\oLD^{\top}$ as follows:
	\begin{gather}
		F\big(C\big)\eq \cC\,,\tag{1}\\[0.4em]
		F\big(\,\twocell{ot}\,\big) \eq \otimes\,,\quad 
		F\big(\,\twocell{oteta}\,\big)\eq 1\,,\quad 
		F\big(\,\twocell{par}\,\big)\eq \partimes \,,\quad 
		F\big(\,\twocell{pareta}\,\big)\eq K \,,\tag{2}\\[0.4em] 
		F\big(a\big)\eq \ao\,,\quad
		F\big(\apc\big)\eq \ap\,,\quad
		F\big(l\big) \eq \lo \,,\quad 
		F\big(\lpc\big) \eq \lp \,,\quad
		F\big(r\big) \eq \ro \,,\quad
		F\big(\rpc\big) \eq \rp \,,\tag{3-1}\\[0.1em]
		F\big(d^l\big)\eq (\distl)^{-1}\,,\quad 
		F\big(d^r\big)\eq \distr\,,\label{3-2}\tag{3-2}\\[0.1em]
		F\big(k^l\big)_{A,B }\eq (\ro_A \prt B) \circ \distl_{A,1,B}\,,\quad 
		F\big(k^r\big)_{A,B} \eq (\rp_A \ot B) \circ (\distr_{A,K,B})^{-1}\,, \quad \text{for } A,B\in \cC\,.\label{3-3}\tag{3-3}
	\end{gather}
	Since $\cC$ is strong, we can use inverse distributors in \eqref{3-2} and \eqref{3-3}.
	
	\smallskip
	
	The assignment $F$ respects the relations in the quotient $\mathbf{LDCat}$: Relations \hyperref[Crit1]{(1)}–\hyperref[Crit4]{(4)} and \hyperref[Crit5]{(5)}–\hyperref[Crit10]{(10)} are satisfied by the commutativity of diagrams \hyperref[eq:T1]{(T1)}–\hyperref[eq:T4]{(T4)} and \eqref{eq:P1}–\eqref{eq:P6}, respectively. Relations \hyperref[Crit11]{(11)}–\hyperref[Crit15]{(15)} follow from Mac Lane's coherence theorem for $\otimes$; \hyperref[Crit16]{(16)}–\hyperref[Crit20]{(20)} from coherence for $\partimes$; and \hyperref[Crit21]{(21)}–\hyperref[Crit22]{(22)} from defining equation \eqref{3-3}. 
	
	The remaining relations \hyperref[Crit23]{(23)}–\hyperref[Crit34]{(34)} are verified using the linearly distributive coherence axioms from Definition \ref{def:LD-cat}, naturality squares, and coherence for $\otimes$ and $\partimes$. Specifically, relation \hyperref[Crit23]{(23)} follows from \hyperref[eq:T4]{(T4)} and naturality of $\rp$; \hyperref[Crit24]{(24)} from \hyperref[eq:T1]{(T1)} and $\lo_1=\ro_1$; \hyperref[Crit25]{(25)} from \eqref{eq:P5}, \hyperref[eq:T2]{(T2)}, and naturality of $\ap$; \hyperref[Crit26]{(26)} from \eqref{eq:P1} and naturality of $\distl$; \hyperref[Crit27]{(27)} from \eqref{eq:P3} and naturality of $\ap$; \hyperref[Crit28]{(28)} from \eqref{eq:P2} and naturality of $\distr$; \hyperref[Crit29]{(29)} from \hyperref[eq:T2]{(T2)} and naturality of $\ro$; \hyperref[Crit30]{(30)} from \hyperref[eq:T3]{(T3)} and $\lp_K=\rp_K$; \hyperref[Crit31]{(31)} from \eqref{eq:P2}, \hyperref[eq:T4]{(T4)}, and naturality of $\ao$; \hyperref[Crit32]{(32)} from \eqref{eq:P6} and naturality of $\distr$; and \hyperref[Crit33]{(33)} from \eqref{eq:P4} and naturality of $\ao$; and \hyperref[Crit34]{(34)} from \eqref{eq:P5} and naturality of $\distl$. 
	
	Relation \hyperref[Crit35]{(35)} follows from the commutativity of the following outer diagram, for $A,B\in \cC$:
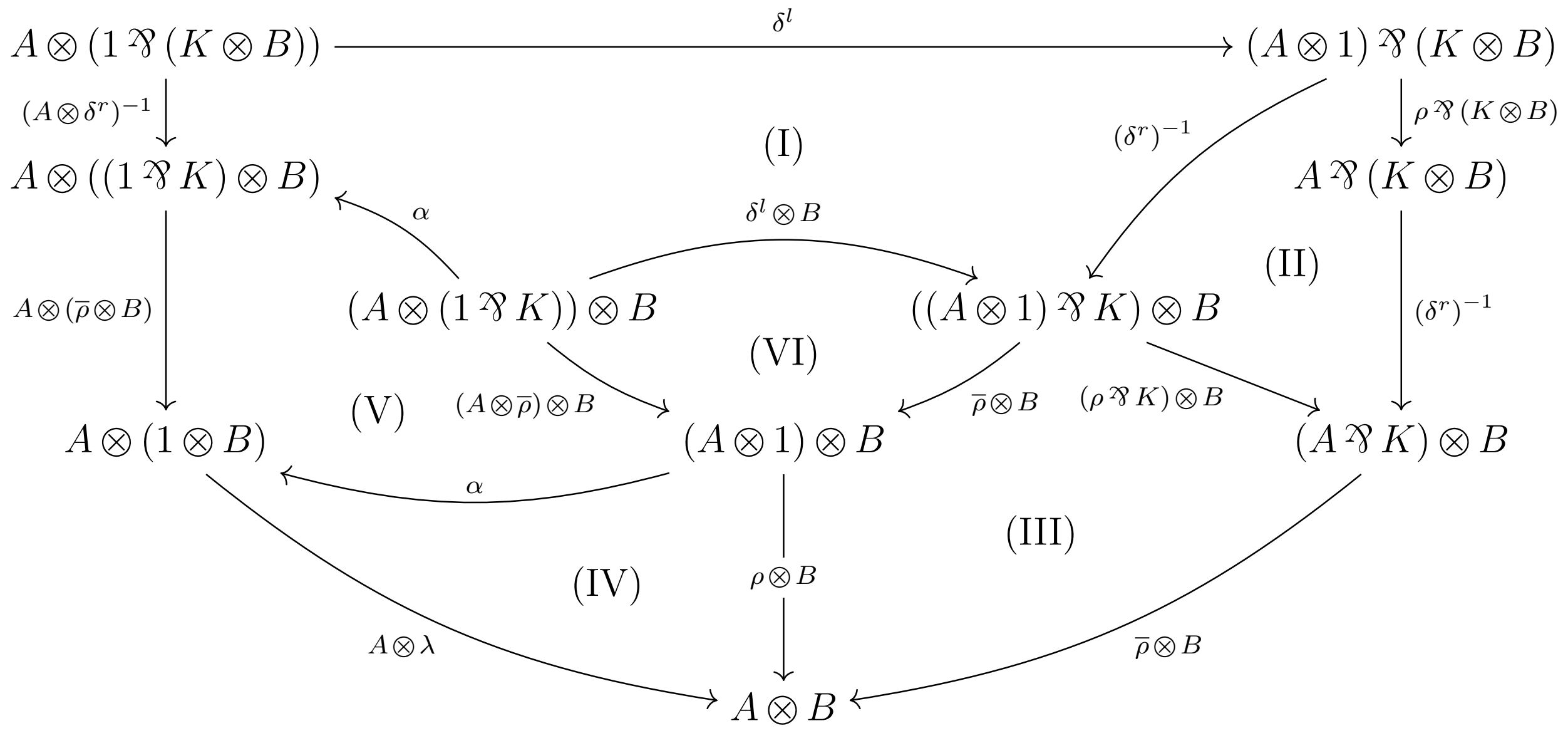
	\begin{equation*}
		\adjustbox{max width=\textwidth}{
			\begin{tikzcd}[cramped,column sep=0.01pt]
				{A\ot (1 \prt {(K \ot B)})} &&&&&& {(A\ot 1) \prt {(K \ot B)}}\\
				{A\ot ((1 \prt K) \ot B)} &&&&&& {A \partimes {(K\otimes B)}} \\ 
				&& {(A\ot (1 \prt K)) \ot B} && {((A\ot 1) \prt K) \ot B} &&\\
				{A\ot (1 \ot B)} &&& {(A\ot 1) \ot B} &&& {(A \partimes K)\otimes B}
				\\
				\\
				\\
				&&& {A \otimes B} &&&
				\arrow["\distl"{yshift=2pt}, bend left=0, from=1-1, to=1-7]
				\arrow["\ro \prt {(K \, \ot \, B)}"{xshift=1pt}, bend left=0, from=1-7, to=2-7]
				\arrow["{(\distr)}^{-1}"{xshift=1pt}, bend left=0, from=2-7, to=4-7]
				\arrow["\rp\, \ot \, B"{xshift=1pt}, bend left=15, from=4-7, to=7-4]
				\arrow["(A\,\ot \,{\distr})^{-1}"'{xshift=-1pt},bend right=0, from=1-1, to=2-1]
				\arrow["A \,\ot \,(\rp\,\ot\,B)"'{xshift=-1pt}, bend right=0, from=2-1, to=4-1]
				\arrow["A\, \ot \, \lo"', bend right=15, from=4-1, to=7-4]
				\arrow["{(\distr)}^{-1}"'{xshift=3pt}, bend right=15, from=1-7, to=3-5]
				\arrow["(\ro \prt K) \,\ot\, B"', from=3-5, to=4-7]
				\arrow["\rp \,\ot \, B", bend left=10, from=3-5, to=4-4]
				\arrow["\ro \, \ot \, B"'{description}, from=4-4, to=7-4]
				\arrow["(A \,\ot\, \rp)\, \ot\, B"', bend right=10, from=3-3, to=4-4]
				\arrow["\distl \,\ot\, B"{yshift=2pt}, bend left=20, from=3-3, to=3-5]
				\arrow["\ao"', bend right=15, from=3-3, to=2-1]
				\arrow["\ao"', bend left=15, from=4-4, to=4-1]
				\arrow[phantom,"\textup{(I)}"{description, yshift=30pt}, from=1-1, to=4-7]
				\arrow[phantom,"\textup{(II)}"{description,yshift=-20pt}, bend left=20, from=1-7, to=3-5]
				\arrow[phantom,"\textup{(III)}"{xshift=20pt}, bend left=20, from=3-5, to=7-4]
				\arrow[phantom,"\textup{(IV)}", bend right=100, from=4-4, to=7-4]
				\arrow[phantom,"\textup{(V)}"{xshift=-30pt,yshift=10pt}, from=2-1, to=7-4]
				\arrow[phantom,"\textup{(VI)}"{yshift=-1pt}, bend left=-10, from=3-3, to=3-5]
			\end{tikzcd}
		}
	\end{equation*}
	Diagram (I) commutes by \eqref{eq:P2}; (II) by naturality of $\distr$; (III) by naturality of $\rp$; (IV) by the triangle identity for $\ot$; (V) by naturality of $\ao$; and (VI) by \hyperref[eq:T4]{(T4)}.  Relation \hyperref[Crit36]{(36)} holds similarly by \eqref{eq:P5}, \hyperref[eq:T2]{(T2)}, the triangle identity for $\prt$, and the naturality of $\ap,\distl,\ro$.
	
	\smallskip
	
	This shows that $F$ is well-defined. Conversely, any $\mathbf{LDCat}$-algebra gives rise to a strong linearly distributive category by relations \hyperref[Crit11]{(11)}, \hyperref[Crit15]{(15)}, \hyperref[Crit16]{(16)}, \hyperref[Crit20]{(20)}, \hyperref[Crit5]{(5)}–\hyperref[Crit10]{(10)}, and \hyperref[Crit1]{(1)}–\hyperref[Crit4]{(4)}.
\end{proof}

We now prove the main result of this note.

\begin{theorem}\label{thm:main theorem}
In a strong linearly distributive category every formal diagram commutes.
\end{theorem}

\begin{proof}
Under the bijection in the proof of Proposition \ref{prop:LDCat-algebras}, formal diagrams in strong linearly distributive categories correspond to images $(F(x_1),F(x_2))$ of $3$-spheres $(x_1,x_2)$ of $\mathbf{LDCat}$. Since $\mathbf{LDCat}$ is aspherical by Propositions~\ref{prop: homotop basis free track} and~\ref{prop: completed computad converg}, the theorem follows.
\end{proof}

In a linearly distributive category, the distributors need not be invertible. Full coherence therefore does not hold, as the following examples show. We begin with an observation.

\begin{remark}\label{rem:adj}
In a linearly distributive category, tensoring with the other unit yields adjoint endofunctors 

\begin{equation*}\label{eq:adjunction-L-R}
	L:= (K \ot -) \;\; \dashv \;\; (1 \partimes -)=:R
\end{equation*}
with counit and unit
\begin{align*}
	LR(A) \sxlongrightarrow{\distl_{K,1,A}} {(K \ot 1)} \partimes A \Simeq A\,,
	\qquad \quad \qquad
	A \Simeq (1 \partimes K)\ot A \sxlongrightarrow{\distr_{1,K,A}} RL(A)\,,
\end{align*}
where the unnamed isomorphisms are unit constraints. The verification of the snake relations for this adjunction relies on the triangle axioms \hyperref[eq:T1]{(T1)}–\hyperref[eq:T4]{(T4)}.
\end{remark}

\begin{example}\label{ex:DRS counterexample}
	Recall from \cite[Ex.~5.9]{DRS26} that, in a linearly distributive category with non-invertible distributors, the two natural transformations involving the counit of the adjunction in Remark~\ref{rem:adj}, depicted in Figure~\ref{fig:capcap}, differ in general.
	
	\begin{figure}[H]
		\centering
		\scalebox{1}{
		\begin{tikzpicture}
			\draw[line width=1] (5.1,0) to[out=90,in=180] (5.5,0.8)
			to[out=0,in=90] (5.9,0);
			\node at (4.85,0) {$L$};
			\node at (5.65,0) {$R$};
			\draw[line width=1] (4.5,0) to (4.5,3);
			\draw[line width=1] (3.5,0) to (3.5,3);
			\node at (3.25,0) {$L$};
			\node at (4.25,0) {$R$};
			\node at (7.4,1.5) {$\Neq$};
			\draw[line width=1] (9,0) to[out=90,in=180] (9.4,0.8)
			to[out=0,in=90] (9.9,0);
			\node at (8.85,0) {$L$};
			\node at (9.65,0) {$R$};
			\draw[line width=1] (11.5,0) to (11.5,3);
			\draw[line width=1] (10.5,0) to (10.5,3);
			\node at (10.25,0) {$L$};
			\node at (11.25,0) {$R$};
		\end{tikzpicture}}
		\caption{Two natural transformations $(L\circ R)^2\to L \circ R$\,.}
		\label{fig:capcap}
	\end{figure}
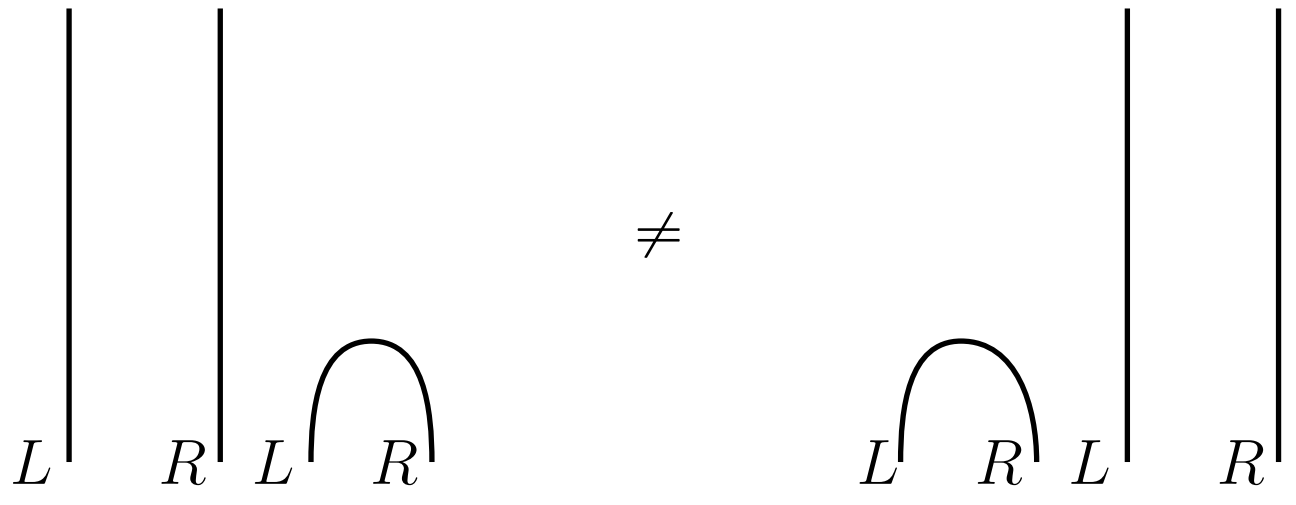
	Postcomposing these two natural transformations with the counit yields the same natural transformation. If the distributors are invertible, then so is the counit, and hence the two natural transformations in Figure~\ref{fig:capcap} agree.
\end{example} 

Another counterexample to coherence in a general linearly distributive category is the following generalization of a variant of the \emph{triple-dual problem} discussed in \cite[Ex.~4.2]{BCST96}. The original triple dual problem appears in \cite[Diag.~(1.4)]{KeMac71}; see also \cite{Hu12,DuVi19}.

\begin{example}\label{ex:counterexample full coherence}
	For positive integers \(m\) and \(n\), the counit and unit from Remark \ref{rem:adj} yield, as shown in Figure~\ref{fig: triple dual problem}, two endomorphisms of the endofunctor \[R\circ (L\circ R)^{m+n-1}.\] These two natural transformations differ in general, e.g. in the linearly distributive category of bimodules over the three-dimensional algebra \(\CC[a,b]/\langle a^2,b^2, ab \rangle\); compare \cite[Ex. 5.9]{DRS26}.
	
	\begin{figure}[H]
		\centering
		\scalebox{1}{
		\begin{tikzpicture}
			\draw[line width=1] (3.1,0) to[out=90,in=180] (3.5,0.8)
			to[out=0,in=90] (3.9,0);
			\node at (2.85,0) {$L$};
			\node at (3.65,0) {$R$};
			\draw[line width=1] (5.1,0) to[out=90,in=180] (5.5,0.8)
			to[out=0,in=90] (5.9,0);
			\node at (4.85,0) {$L$};
			\node at (5.65,0) {$R$};
			\fill (4.3,0.4) circle (1pt);
			\fill (4.45,0.4) circle (1pt);
			\fill (4.6,0.4) circle (1pt);
			\draw[line width=1] (0,0) to[out=90, in=270] (3.5,3);
			\node at (-0.25,0) {$R$};
			\node at (3.25,3) {$R$};
			\draw[line width=1] (0.8,0) to[out=90, in=200] (2.8,1.5) to[out=20, in=270] (4.3,3);
			\node at (0.55,0) {$L$};
			\node at (4.05,3) {$L$};
			\draw[line width=1] (2.3,0) to[out=90, in=270] (6.1,3);
			\node at (5.85,3) {$R$};
			\node at (2.05,0) {$R$};
			\fill (3.3,1.5) circle (1pt);
			\fill (3.45,1.5) circle (1pt);
			\fill (3.6,1.5) circle (1pt);
			\draw[line width=1] (0,3) to[out=270,in=180] (0.4,2.2) to[out=0,in=270] (0.8,3);
			\node at (-0.25,3) {$R$};
			\node at (0.55,3) {$L$};
			\fill (1.2,2.6) circle (1pt);
			\fill (1.35,2.6) circle (1pt);
			\fill (1.5,2.6) circle (1pt);
			\draw[line width=1] (2,3) to[out=270,in=180] (2.4,2.2)
			to[out=0,in=270] (2.8,3);
			\node at (1.75,3) {$R$};
			\node at (2.55,3) {$L$};
			\draw[line width=1, decorate, decoration={brace, mirror,amplitude=8pt}]
			(0,-0.3) -- (2.3,-0.3);
			\draw[line width=1, decorate, decoration={brace, mirror,amplitude=8pt}]
			(3.1,-0.3) -- (5.9,-0.3);
			\draw[line width=1, decorate, decoration={brace,amplitude=8pt}]
			(0,3.3) -- (2.8,3.3);
			\draw[line width=1, decorate, decoration={brace,amplitude=8pt}]
			(3.5,3.3) -- (6.1,3.3);
			\node[font=\footnotesize] at (4.6,-1.15) {$m$ counits};
			\node[align=center, font=\footnotesize] at (1.1,-1.15) {$2n-1$\\ through-strands};
			\node[font=\footnotesize] at (1.5,4.05) {$m$ units};
			\node[align=center, font=\footnotesize] at (4.8,4.1) {$2n-1$\\ through-strands};
			\node at (7.4,1.5) {$\Neq$};
			\draw[line width=1] (9.1,0) to (9.1,3);
			\node at (8.75,0) {$R$};
			\draw[line width=1] (9.9,0) to (9.9,3);
			\node at (9.65,0) {$L$};
			\fill (12.1,1.5) circle (1pt);
			\fill (12.25,1.5) circle (1pt);
			\fill (12.4,1.5) circle (1pt);
			\draw[line width=1] (15.2,0) to (15.2,3);
			\node at (14.95,0) {$R$};
			\draw[line width=1, decorate, decoration={brace, mirror,amplitude=8pt}]
			(9,-0.3) -- (15.3,-0.3);
			\draw[line width=1, decorate, decoration={brace,amplitude=8pt}]
			(9,3.3) -- (15.3,3.3);
			\node[align=center, font=\footnotesize] at (12.1,-1.15) {$2m+2n-1$\\ through-strands};
			\node[align=center, font=\footnotesize] at (12.1,4.2) {$2m+2n-1$\\ through-strands};
		\end{tikzpicture}}
		\caption{Two endomorphisms of the endofunctor $R\circ (L\circ R)^{m+n-1}$.}
		\label{fig: triple dual problem}
	\end{figure}
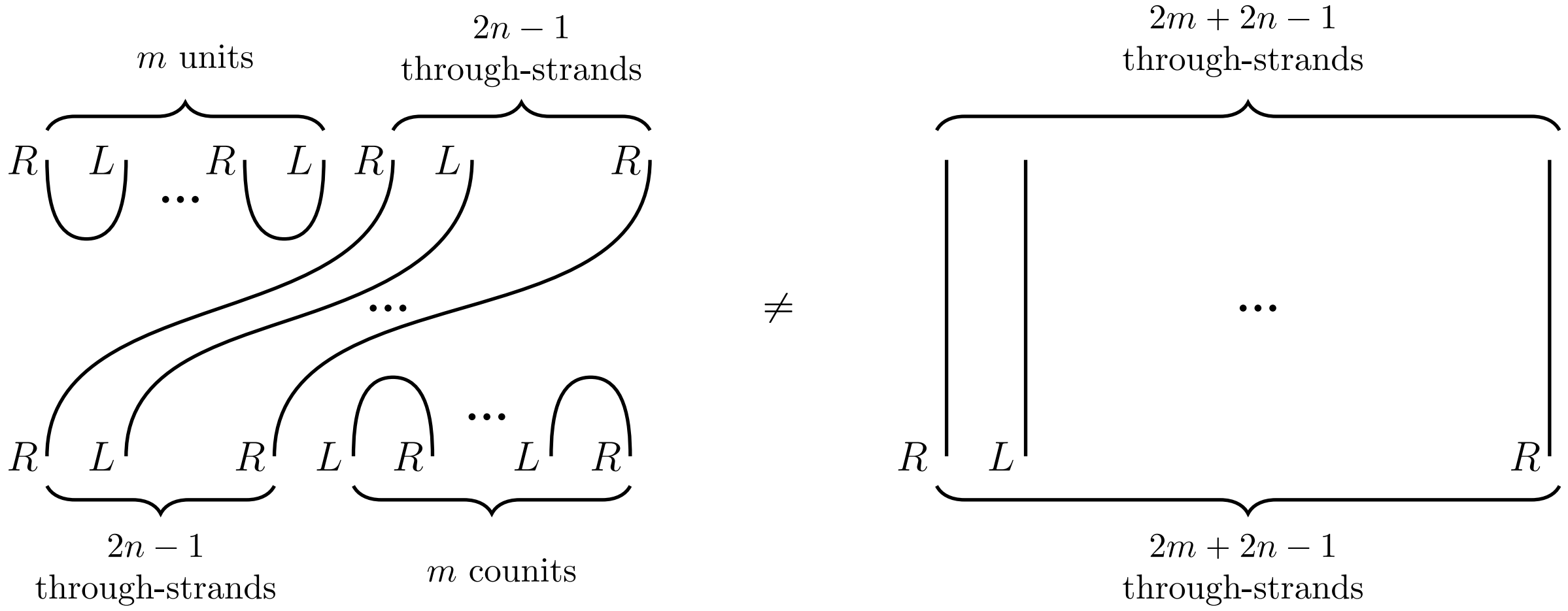
	Denote the left-hand natural transformation in Figure \ref{fig: triple dual problem} by $\nu(m,n)$. For positive integers $m,n,k,l$ with $m+n=k+l$, one finds	
\begin{equation*}
{\nu(k,l)} \circ {\nu(m,n)} \eq \nu(r,s)\,,
\end{equation*}
 where 
\begin{equation*}
	r\eq m+\operatorname{min}(k,n-1) \qquad \text{and} \qquad s \eq \operatorname{max}(1,n-k)\,.
\end{equation*}
 This follows by distinguishing the cases $n-k\leq 1$ and $n-k > 1$, and using the snake relations for the adjunction \hyperref[eq:adjunction-L-R]{$L \dashv R$}. In particular, setting $k=m+n-1$, gives 
 \begin{equation*}
 	\nu(k,1)\circ\nu(m,n)\eq \nu(k,1) \,.
 \end{equation*}
 If the distributors are isomorphisms, then the idempotent $\nu(k,1)$ is invertible. Therefore, ${\nu(m,n)=\id}$ for all positive integers $m,n$, so the two endomorphisms in Figure~\ref{fig: triple dual problem} agree.
\end{example}

\begin{center}
	\sc AI disclosure\\[.5em]
	\normalfont No AI was used for research. AI was used sparingly to suggest alternative wording.
\end{center}

\medskip

\begin{center}
	\sc Acknowledgments\\[.5em]
\end{center}
All string diagrams were typeset with Yves Guiraud’s Catex package.

M.D. and C.S. acknowledge support from the DFG through the CRC 1624 \emph{Higher Structures, Moduli Spaces and Integrability}, project number 506632645. C.S. is partially funded by the Excellence Cluster EXC 2121 \emph{Quantum Universe}, project number 390833306.

\bibliographystyle{alphaurl}
\bibliography{references}
\end{document}